\documentclass[a4paper,leqno]{article}
\usepackage{amsmath,amsthm,amssymb,amscd}

\chardef\bslash=`\\ 

\theoremstyle{change}
\newtheorem{thm}{Theorem}[section]
\newtheorem{cor}[thm]{Corollary}
\newtheorem{lem}[thm]{Lemma}
\newtheorem{prop}[thm]{Proposition}

\theoremstyle{definition}
\newtheorem{defn}{Definition}[section]

\theoremstyle{remark}
\newtheorem{rem}{Remark}[section]
\newtheorem*{notation}{Notation}

\numberwithin{equation}{section}

\newcommand{\maprU}[1]{
\smash{\mathop{%
\hbox to 1cm{\rightarrowfill}}\limits^{#1}}}
\newcommand{\maprD}[1]{
\smash{\mathop{%
\hbox to 1cm{\rightarrowfill}}\limits_{#1}}}
\newcommand{\maplU}[1]{
\smash{\mathop{%
\hbox to 1cm{\leftarrowfill}}\limits^{#1}}}
\newcommand{\maplD}[1]{
\smash{\mathop{%
\hbox to 1cm{\leftarrowfill}}\limits_{#1}}}

\newcommand{\mapdL}[1]{
\Big\downarrow\llap{$\vcenter{\hbox{$\scriptstyle#1\,\,\,$}}$}}
\newcommand{\mapdR}[1]{
\Big\downarrow\rlap{$\vcenter{\hbox{$\scriptstyle#1$}}$}}

\newcommand{\leseq}{\mbox{\underline{$<$}}}

\newcommand{\eval}[2][\right]{\relax
  \ifx#1\right\relax \left.\fi#2#1\rvert}

\begin{document}
\title{Rational integrals of the second kind \\
on a complex projective manifold \\
and its primitive cohomology
\thanks{This work is supported by the Grant-in-Aid for Scientific Research (No.
19540093), The Ministry of Education, Science and Culture, Japan} }
\author{Shoji TSUBOI \\
{\small Department of Mathematics and Computer Science, Kagoshima University}\\
{\footnotesize e-mail: tsuboi@sci.kagoshima-u.ac.jp}}
\date{}
\maketitle

\noindent
{\sc Abstract}:
Let $X$ be a complex algebraic manifold of dimension $n+1$ 
 embedded in a sufficiently higher dimensional complex projective space $\Bbb{P}^{N}(\Bbb{C})$, and $Y$ a generic hyperplane section of $X$. By sheaf cohomological method, we prove the well-known facts that the primitive cohomology group $H^p(X,\mathbb{C})_0$\, ($1\leseq p\leseq n+1$) is isomorphic to the De Rham cohomology group $I^p(X,(p+1)Y)_0$ of closed rational $p$-forms of the 2nd kind on $X$, having poles of order $p+1$ (at most) along $Y$ only, and that the {\it Hodge filtration} of $H^p(X,\mathbb{C})_0$\ is isomorphic to the one of $I^p(X,(p+1)Y)_0$ defined by the order of poles along $Y$. On the other hand, we have a long exact sequence of cohomology\par
\vskip 7pt

$\to H^p(X,\mathbb{C})\xrightarrow{r^p} H^p(X-Y,\mathbb{C})\xrightarrow{R^p} 
H^{p-1}(Y,\mathbb{C})\xrightarrow{G^{p-1}} H^{p+1}(X,\mathbb{C})\to \cdots$,
\par
\vskip 7pt 
\noindent
which is dual to \par
\vskip 7pt

$\to H_p(X,\mathbb{C})\xleftarrow{\iota_p} H^c_p(X-Y,\mathbb{C})
\xleftarrow{\tau_{p-1}}
H_{p-1}(Y,\mathbb{C})\xleftarrow{G_{p+1}} H_{p+1}(X,\mathbb{C})\to \cdots$, 
 \par
\vskip 7pt
\noindent 
where $H^c_*$ denotes compact support homology group (cf. (\ref{e1.2})). 
Using these exact sequences, we describe the mixed Hodge structure on $H^p(X-Y,\mathbb{C})$ and the Hodge filtration of the middle primitive cohomology group $H^{n}(Y,\mathbb{C})_0$ of $Y$ in terms of rational integrals on $X$.
\par

\bigskip
\bigskip
\tableofcontents
\bigskip
\bigskip
\bigskip
\section*{Summary}
Let $X$ be a non-singular irreducible algebraic variety of dimension $n+1$ 
 embedded in a sufficiently higher dimensional complex projective space $\Bbb{P}^{N}(\Bbb{C})$, and $Y$ a generic hyperplane section of $X$. We shall use the following notation:
\vskip 5pt

\noindent
$\Omega_X^q$ : the sheaf of germs of holomorphic $q$-forms on $X$,\par
\vskip 2pt

\noindent
$\Omega_X^q(kY)$ : the sheaf of germs of meromorphic $q$-forms having poles of order $k$ \\
\hskip 20pt (at most) along $Y$ as their only singularities on $X$,\par
\vskip 2pt

\noindent
$\Omega_X^q(*Y)$ : the sheaf of germs of meromorphic $q$-forms having poles of arbitrary \\
\hskip 20pt order along $Y$ as their only singularities on $X$,\par
\vskip 2pt

\noindent
$\Omega_X^q(\log Y)$ : the sheaf of germs of meromorphic $q$-forms having logarithmic \\
\hskip 20pt poles (at most) along $Y$ as their only singularities on $X$.\par
\vskip 5pt

\noindent
We denote by $\Phi^q_X$, $\Phi_X^q(kY)$, e.t.c., the subsheaves consisting of closed forms of each ones.  
On the complex $\Omega_X^\cdot$ we define a decreasing filtration $F=\{F^k\}_{0\leseq k\leseq n+1}$ (the Hodge filtration)  by 
the subcomplexes
\begin{equation*} 
F^k(\Omega_X^\cdot)^q=
\left\{
 \begin{array}{cl}
  0 &  q<k \\
  \Omega_X^q &  k\leseq q. 
  \end{array}
 \right.
\end{equation*}

\noindent
On the complex $\Omega_X^\cdot(\log Y)$ we define the Hogde filtartion similarly, 
and another increasing filtration $W=\{W_0\subset W_1\}$ (the weight filtration) by 
\begin{equation*}
 W_0(\Omega_X^\cdot(\log Y))=\Omega_X^\cdot, \quad 
W_1(\Omega_X^\cdot(\log Y))=\Omega_X^\cdot(\log Y).
\end{equation*}

\noindent
Then $(\Omega_X^\cdot, F)$ becomes the {\it cohomological Hodge complex}, and 
$(\Omega_X^\cdot(\log Y), W, F)$ the {\it cohomological mixed Hodge complex} (cf. \S 3). 
They induce the Hodge structure on the cohomology $H^p(X,\mathbb{C})$, and the mixed Hodge structure on the cohomology $H^p(X-Y,\mathbb{C})$.  We define 
\begin{equation*}
 I^p_k(X,(p+1)Y):=\frac{\Gamma(X,\Phi_X^p((p-k+1))Y)}{d\Gamma(X,\Omega_X^{p-1}((p-k))Y)} \quad (0\leseq k\leseq p)
 \end{equation*}

\noindent
and denote by $I_k^p(X,(p+1)Y)_0$ the subspace of $I_k^p(X,(p+1)Y)$ generated 
by closed moromorphic $p$-forms of the {\it second kind\/} (cf. Definition \ref{def2.2}).  Assume that
\begin{equation*}
 H^p(X,\Omega_X^q(kY))=0 \quad \text{for} \quad p\geq 1, q\geq 0 \quad \text{and} \quad k\geq 1.
 \end{equation*}
 \noindent
 Then we have 
\begin{eqnarray*}
& & F^kH^p(X-Y,\mathbb{C})\simeq I^p_k(X,(p+1)Y) \quad 0\leseq k\leseq p, \\
& & F^kH^p(X,\mathbb{C})_0 \simeq I^p_k(X,(p+1)Y)_0 \quad 0\leseq k\leseq p, \\
& & Gr_q^{W[q]}H^q(X-Y,\mathbb{C})=W[q]_qH^q(X-Y,\mathbb{C})=I^q(X,\ast Y)_0,\\
& & Gr_{q^+1}^{W[q]}H^q(X-Y,\mathbb{C})=I^q(X,\ast Y)/I^q(X,\ast Y)_0, \\
& & F^kGr_{q}^{W[q]}H^q(X-Y,\mathbb{C})\simeq F^kH^q(X,\mathbb{C})_0, \quad \text{and}\\
& & F^kGr_{q+1}^{W[q]}H^q(X-Y,\mathbb{C})
\simeq{\rm Ker}\{F[-1]^kH^{q-1}(Y,\mathbb{C})_0\xrightarrow{G}F^kH^{q+2}(Y,\mathbb{C})\},
\end{eqnarray*}
\smallskip

\noindent
where $H^p(X,\mathbb{C})_0$ denotes the $p$-th primitive coholomology of $X$,
$F^k$ the $k$-th Hodge filtration of cohomology, and $W[q]$ the shift to the right on the degree of $W$ by $q$.  (Theorem \ref{thm3.1}, Theorem \ref{thm3.3} and Proposition \ref{prop2.3}). Furthermore, let $Y^\prime$ be a generic hypersurface of $\mathbb{P}^{N}(\mathbb{C})$ of sufficiently higher degree so that 
\begin{equation*}
 H^p(Y,\Omega_Y^q(kZ))=0 \quad \text{for} \quad p\geq 1, q\geq 0 \quad \text{and} \quad k\geq 1,
\end{equation*}

\noindent
where $Z=Y\cdot Y^\prime$.  Then we can define the {\it generalized 
Poincar{\'e} r{\'e}sidue map}
\begin{equation*}
R\acute{e}s: I^{n+1}(X,(n+2)Y)\to I^{n}(Y,(n+1)Z)_0
\end{equation*}

\noindent
and prove that 
\begin{eqnarray*}
F^kH^n(Y,\mathbb{C})_0&\simeq& I^n_k(Y,(n+1)Z)_0 \\
                      &\simeq& R\acute{e}s(I^{n+1}_{k+1}(X,(n+2)Y))\oplus 
                      r^n(I^n_k(X,(n+1)Y^\prime)_0)),
\end{eqnarray*}

\noindent
where $r^n$ denotes the map induced by the natural map $H^n(X,\mathbb{C})_0\to H^n(Y,\mathbb{C})_0$ (Thorem \ref{thm4.1}). These results might be considered as a generalization of those by P. A. Griffith in the case of a hypersurface in a complex projective space (cf. \cite{Gri-1}).\\
\bigskip
\bigskip

\section{Some remarks on primitive cohomology and homology 
of algebraic manifolds}

 Let {X} be a non-singular irreducible algebraic variety of dimension $n+1$ 
 embedded in a higher dimensional complex projective space $\Bbb{P}^{N}(\Bbb{C})$ and $Y$ a generic hyperplane section of $X$. In what follows we call such $Y$ a {\it prime section} of $X$.  We denote by $\Omega$ 
 the restriction to $X$ of the fundamental form of the Fubini-Study metric on 
$\Bbb{P}^{N}(\Bbb{C})$. $\Omega$ is a closed 2-form whose cohomology class $[\Omega]\in H^2(X,\Bbb{C})$ is the Poincar\'e dual of the homology class $[Y]\in H_{2n}(X,\Bbb{C})$ associated to the the prime section $Y$. We define $L(\omega):=\Omega\wedge\omega$ for a ($\Bbb C$-valued) $C^\infty$ diferential $q$-form $\omega$ on $X$. If $\omega$ is a closed form (resp. detived form), then $L(\omega)$ is also a closed form (resp. derived form) for $\Omega$ is a closed form.  Hence $L$ define a homomorphism $H^q(X,\Bbb{C})\to H^{q+2}(X,\Bbb{C})$\,\,($0\leseq q\leseq 2n$). Throughput this paper we always idetify the ordinary cohomology with 
the De Rham cohomology. We call this cohomology operator {\it Hodge operator\/} and denote it by the same letter $L$. \par
\smallskip

\begin{defn} \label{d.1.1} A\, $C^{\infty}$ differential $q$-form \,\,($0\leseq q\leseq n+1$) $\omega$ is said to be {\it primitive\/} if $L^{n-q+2}(\omega)=0$ \,\
($L^{n-q+2}=\underbrace{L\circ\cdots\circ L}_{\mbox{$\scriptsize{n-q+2}$ \scriptsize{times}}}$). A (De Rham) cohomology class containing a closed, primitive $C^{\infty}$ differential form is said to be a {\it primitive\/} cohomology class.  
\end{defn}

 We call the subgroup of $H^q(X,\Bbb{C})$ which consists of all primitive cohomology classes the {\it $q$-th primitive cohomology group\/} of $X$, which we denote by $H^q(X,\Bbb{C})_0$.
 
\begin{rem} Originarlly, a \,$C^{\infty}$ differetial $q$-form ($\leseq q\leseq n+1$) $\omega$ on $X$ is defined to be primitive if $\Lambda\omega=0$, $\Lambda$ is the adjoint operator of $L$ with respect to the Hodge metric on $X$ which is the restriction of the Fubini-Study metric on $\Bbb{P}^{N}(\Bbb{C})$.  The above definition of primitive forms is equivalent to the original one (cf. \cite{G-H}).
\end{rem}

	The following facts are fundamental for the Hodge operator $L$. 

\begin{thm}{\rm (Hard Lefshets Theorem)}\label{thm1.1}
$$
   L^k: H^{n+1-k}(X,\Bbb{C})\simeq  H^{n+1+k}(X,\Bbb{C})\quad (1\leseq k\leseq n+1)
$$
\end{thm}

\begin{thm}{(\rm Lefshets decomposition)}\label{thm1.2}
\begin{itemize}
\item[{\rm (i)}] $L: H^{q-2}(X,\Bbb{C})\to H^{q}(X,\Bbb{C})$ is injective and \par
        $H^{q}(X,\Bbb{C})\simeq LH^{q-2}(X,\Bbb{C})\oplus H^{q}(X,\Bbb{C})_0$ 
        \quad ($2\leseq q\leseq n+1)$. 
\item[{\rm (ii)}] $H^{n+1+k}(X,\Bbb{C})\simeq L^kH^{n+1+k}(X,\Bbb{C})_0\oplus L^{k+1}H^{n-1-k}(X,\Bbb{C})$
\end{itemize}
\end{thm}

	By restriction\, $C^\infty$ differential $q$-forms on $X$ to $Y$, we obtain a cohomology map $r^q: H^{q}(X,\Bbb{C})\to H^{q}(Y,\Bbb{C})$, for which the folowing holds.

\begin{thm} {\rm (Weak Lefshetz Theorem)} \label{thm1.3}
\begin{itemize}
\item[{\rm (i)}] $r^q: H^{q}(X,\Bbb{C})\simeq H^{q}(Y,\Bbb{C})$
        \quad ($0\leseq q\leseq n-1)$. 
\item[{\rm (ii)}] $r^n: H^{n}(X,\Bbb{C})\to H^{n}(Y,\Bbb{C})$ is 
injective.
\end{itemize}
\end{thm}

\noindent
For the proofs of the theorems above we refer to \cite{G-H}. 

\begin{cor}\label{cor1.4}
$$
0\to H^{n+1}(X,\Bbb{C})_0\to H^{n+1}(X,\Bbb{C})\xrightarrow{r^{n+1}} H^{n+1}(Y,\Bbb{C})\to 0. \quad \mbox{(exact)}
$$
\end{cor}

\begin{proof} By (\ref{thm1.2}), (i) and (\ref{thm1.1}), we have 
\begin{eqnarray*}
 \begin{array}{ccccc}
  0 & \to & H^{n-1}(X,\Bbb{C})& \xrightarrow{L} & H^{n+1}(X,\Bbb{C}) \\
    &     & \mapdL{\simeq} &  &  \mapdR{r^{n+1}} \\
    &     &  H^{n-1}(Y,\Bbb{C}) & 
\xrightarrow[\simeq]{L} 
    & H^{n+1}(Y,\Bbb{C})
 \end{array}
 \end{eqnarray*}
 \offinterlineskip
 \noindent
 and, 
 \begin{eqnarray*}
 H^{n+1}(X,\Bbb{C})=H^{n+1}(X,\Bbb{C})_0\oplus LH^{n-1}(X,\Bbb{C}). 
 \end{eqnarray*}
 
 \offinterlineskip
 \noindent
 Therefore, 
 \begin{eqnarray*}
 Ker~r^{n+1}=H^{n+1}(X,\Bbb{C})_0
\end{eqnarray*}
\end{proof} 

\begin{cor}\label{cor1.5}
\begin{eqnarray*}
 0\to H^{n+1}(X,\Bbb{C})_0\to H^{n+1}(X,\Bbb{C})\xrightarrow{r^{n+1}}H^{n+1}(Y,\Bbb{C})\to 0\quad \mbox{{\rm (exact)}}
 \end{eqnarray*}
\end{cor}
\bigskip

	In what follows, homology and cohomology are with coefficient in the complex number field if otherwise explicitly mentioned.  Taking a topological tublar neighborhood $U$ of $Y$ in $X$, we consider the homology exact sequence concerning a pair of the topological spaces $(X,X-U)$, 
which is written as follows:
\begin{equation}\label{e1.1}
\begin{split}
&\cdots \to H_q^c(X-U)\xrightarrow{i_q} H_q(X) \\
 &\hskip 50pt \xrightarrow{j_q}H_q(X,X-U)
\xrightarrow{\partial_q}H_{q-1}^c(X-U)\to \cdots,
\end{split}
\end{equation}

\noindent
where $H_\ast^c$ denotes compact support homology groups. Since $X-U$ is a 
deformation retract of $X-U$, $H_{q}^c(X-U)\simeq H_{q}^c(X-Y)$.  By the 
{\it excision axiom\/}, $H_{q}^c(X,X-U)\simeq H_{q}^c(U,\partial U)$.  By 
the {\it Thom isomorphism\/}, $H_{q}^c(U,\partial U)\simeq H_{q-2}^c(Y)$ for 
$q\geq 2$.  We obviously have $H_{q}(U,X-U)=0$ for $0\leseq q\leseq 1$.  
Therefore the homology exact sequence (1.1) is rewritten as follows:
\begin{equation}\label{e1.2}
\begin{split}
&\cdots \to H_q^c(X-Y)\xrightarrow{\iota_q} H_q(X) \\
 &\hskip 50pt \xrightarrow{G_{q}}H_{q-2}(Y)
\xrightarrow{\tau_{q-2}}H_{q-1}^c(X-Y)\to \cdots,
\end{split}
\end{equation}
\noindent
where 
\begin{description}
\item[(i)] the map $\iota_q: H_q^c(X-U)\rightarrow H_q(X)$ is the one induced by the natural inclusion map $\iota: X-Y\to X$, 

\item[(ii)] the map $G_{q}: H_q(X)\rightarrow H_{q-2}(Y)$ is the one which assignes each q-cycle on $M$ to its intersection cycle with $Y$, and 

\item[(iii)] the map ${\tau_{q-2}}: H_{q-2}(Y)\rightarrow H_{q-1}^c(X-U)$ is the one which assighns each $(q-2)$ cycle on $Y$, say $\gamma$, to the cycle 
$\partial U_{\vert\gamma}$ on $X-Y$, the restriction of $\partial U$ over 
$\gamma$.

\end{description}

 In the subsequence we denote the cycle $\partial U_{\vert\gamma}$ in (iii) above by $\tau(\gamma)$.  Taking the cohomology exact sequence dual to (\ref{e1.2}), we have
\begin{equation}\label{e1.3}
\begin{split}
&\cdots \to H^q(X-Y)\xleftarrow{r^q} H^q(X) \\
 &\hskip 50pt \xleftarrow{G^{q-2}}H^{q-2}(Y)
\xleftarrow{R^{q-1}}H^{q-1}(X-Y)\to \cdots,
\end{split}
\end{equation}
 
\noindent
Here the map $G^{q-2}: H^{q-2}(Y)\to H^q(X)$ is the so-called {\it Gysin map\/}.  We are now going to describe the Gysin map $G^{q-2}$ by use of differential forms. We take a sufficiently fine, finite open covering 
${\cal U}=\{U_i\}_{i\in I}$ of $X$ such that, in each open subset $U_i$, $Y$ is defined by a holomorphic equation $\sigma_i=0$.  We put $t_{ij}=\sigma_i/\sigma_j$ for each pair of indexes $(i,j)$ with $U_i\cap U_j\not=\emptyset$.  Then the system of transition functions, with respect to the covering ${\cal U}$, of the line bundle $[Y]$ associated to $Y$ are given by $\{t_{ij}\}$, and 
$\sigma=\{\sigma_i\}$ give rise to a cross-section of $[Y]$ whose zero locus is $Y$.  We take a system $\{a_i\}$ of real positive functions $a_i$ of class $C^\infty$ defined in $U_i$, respectively, satisfying 
\begin{equation*}
  \frac{a_i}{a_j} = \vert t_{ij}\vert^2, \quad \mbox{in}\,\,U_i\cap U_j\not=\emptyset.
\end{equation*}

\noindent
The system $\{a_i\}$ defines a fiber metric on the line bundle $[Y]$.  The 
length function $\vert\sigma\vert$ of the cross-section $\sigma=\{\sigma_i\}$ of $[Y]$ with respect to this fiber metric is given by 
\begin{eqnarray*}
\vert\sigma\vert &=& \sqrt{\sigma_i a_i\overline{\sigma_i}}\\
                 &=& \vert\sigma_i\vert\sqrt{a_i}
\end{eqnarray*}

\noindent
in each $U_i$.  Note that $\vert\sigma\vert^2$ is a globally defined real 
non-negative function of class $C^\infty$.  We define 
\begin{eqnarray*}
\eta &:=& \frac{1}{2\pi i}\partial \log \vert\sigma\vert^2, \\
 \omega &=& \overline{\partial}\eta = 
 \frac{1}{2\pi i}\overline{\partial}\partial \log \vert\sigma\vert^2.
\end{eqnarray*}

\noindent
On each $U_i$, $\eta$ and $\omega$ are written as 
\begin{eqnarray*}
\eta &:=& \frac{1}{2\pi i}(d \log \sigma_i + \partial\log a_i), \\
 \omega &=& \frac{1}{2\pi i}\overline{\partial}\partial \log a_i.
\end{eqnarray*}

\noindent
Note that $\omega$ is a globally defined closed $C^\infty$ form of type $(1,1)$ on $X$, representing the first Chern class $c_1([Y])$ of the line bundle $[Y]$. We denote by $A^\ast(X)$, $A^\ast(X-Y)$ and $A^\ast(Y)$ the De Rham complexes of $\Bbb C$-valued, $C^\infty$ differential forms on $X$, $X-Y$ and $Y$, respectively. 

\begin{defn}\label{d1.2}
$A^\ast(\log Y)$ is defined to be the sub-complex of $A^\ast(X-Y)$ generated 
by $A^\ast(X)$ and $\eta$. 
\end{defn}

 A form $\varphi\in A^\ast(\log Y)$ may be (non-uniquely) written as 
 
\begin{equation}\label{e1.4}
 \varphi=\alpha\wedge \eta + \beta
\end{equation}
 
\noindent
where $\alpha, \beta \in A^\ast(V)$.  The restriction $\alpha_{\vert Y}\in A^\ast(Y)$ is, however, not anbiguous.  Hence we may define $R^\ast: A^\ast(\log Y)\to A^{\ast-1}(Y)$ by 
 \begin{equation}\label{e1.5}
 R^\ast(\varphi):=2\pi\sqrt{-1}\alpha_{\vert Y},
 \end{equation}
 
\noindent
which we call {\it R\'esidue map\/}.  Let $W^\ast\subset A^\ast(\log Y)$ be 
the kernel of $R^\ast$.  There is an obvious inclusion 
\begin{equation*}
A^\ast(X)\mbox{\vbox{\hbox{$\iota$}\hbox{$\subset$}}}W^\ast
\end{equation*}
\smallskip

\begin{prop}\label{prop1.6}
The inclusion $\iota$ induces isomorphisms on $d$ and $\overline{\partial}$ 
cohomologys.
\end{prop}

For the proof we refer to (\cite{Gri-1}), p.49$\sim$p.50. 
\smallskip

\begin{prop}\label{prop1.7}
The Gysin map $G^{q-2}:H^{q-2}(Y,\mathbb{C})\to H^{q}(X,\mathbb{C})$ is described 
using differential forms as follows: For $\alpha\in A^{q-2}(Y)$, choose 
$\tilde{\alpha}\in A^{q-2}(X)$ with $\tilde{\alpha}_{\vert Y}=\alpha$ and set
\begin{equation*}
\gamma(\alpha)=d(\tilde{\alpha}\wedge\eta)=d\tilde{\alpha}\wedge\eta\wedge\eta
+(-1)^{q-2}\tilde{\alpha}\wedge\omega.
\end{equation*}

\noindent
If $\alpha$ is a closed form (resp. deived from), then $\gamma(\alpha)$ is a 
closed form (resp. derived form) in $W^q$.   Furthermore $\gamma(\alpha)$ is 
independent of the choice of $\tilde{\alpha}$ modulo derived form in $W^q$.  Hence, by virtue of Proposition \ref{prop1.7}, the correspondence $[\alpha]\to [\gamma(\alpha)]$ defines a map 
\begin{equation*}
  H^{q-2}(Y,\mathbb{C})\simeq H^{q-2}(A^\ast(X))\to 
  H^{q}(X,\mathbb{C})\simeq H^q(W^\ast),
\end{equation*}

\noindent
which coincides, up to a factor of $\pm 1$, with the Gysin map $G$. 
\end{prop}

\begin{proof}
By the definition of $W^ast$, $\gamma(\alpha)\in W^ast$.  It is obvious that 
if $\alpha$ is a closed form, then $\gamma(\alpha)$ is also closed in $W^q$. 
Assume $\alpha$ is wriiten as $d\beta=\alpha$ for $\beta\in A^{q-3}(Y)$.  
We choose $\tilde{\beta}\in A^{q-3}(X)$ with $\tilde{\beta}_{\vert Y}=\beta$ 
and set 
\begin{equation*}
\xi=(\tilde{\alpha}-d\tilde{\beta})\wedge \eta + (-1)^{q-2}\tilde{\beta}\wedge d\eta.
\end{equation*}

\noindent
Then $\xi\in W^{q-1}$ and 
\begin{eqnarray*}
d\xi &=& d\tilde{\alpha}\wedge\eta + (-1)^{q-2}(\tilde{\alpha}-d\tilde{\beta})\wedge d\eta + (-1)^{q-2}d\tilde{\beta}\wedge d\eta \\
 &=& d\tilde{\alpha}\wedge\eta + (-1)^{q-2}\tilde{\alpha}\wedge d\eta\\
 &=& \gamma(\alpha)
\end{eqnarray*}

\noindent
Thus $\gamma(\alpha)$ is a derived form in $W^\ast$. \par
\smallskip

 The fact that $\gamma(\alpha)$ is independent of the choice of $\tilde{\alpha}$ modulo derived forms in $W^\ast$ is almost trivial.  In fact, if $\tilde{\alpha}'$ is another form in $A^{q-2}(X)$ with $\tilde{\alpha}'_{\vert Y}=\alpha$, then 
 $(\tilde{\alpha}-\tilde{\alpha}')\wedge \eta\in W^{q-1}(X)$ and $d((\tilde{\alpha}-\tilde{\alpha}')\wedge \eta)=d\tilde{\alpha}\wedge \eta-d\tilde{\alpha}'\wedge \eta$, which shows $\gamma(\alpha)$ is uniquely determined up to derivede forms in $W^\ast$.  we wre now going to show that the correspondence $[\alpha]\to 
 [\gamma(\alpha)]$ coincides with the Gysin map $G$.  To do this it suufices to show that for any $q$-cycle $c_q$ on $X$, the integral $\int_{\Gamma}\gamma(\alpha)$ converges and 
\begin{equation}\label{e1.6}
\int_{c_q}\gamma(\alpha)=\pm\int_{c_q\cdot Y}\alpha
\end{equation}

\noindent
holds, where $\Gamma\cdot Y$ denotes the intersection cycle of $\Gamma$ with $Y$. We may assume that $c_q$ intersects $Y$ normally in a $(q-2)$ cycle $c_{q-2}$ with respect to some given hermitian metric on $X$.  For a sufficiently small positive $\varepsilon$, we take a {\it $\varepsilon$-tube} with axis $c_{q-2}$, and lying in $c_q$, normally, 
\begin{equation*}
T_{\varepsilon}(c_{q-2}):=\{~p\in c_q~ \vert ~d_X(p,c_{q-2})\leseq \varepsilon~\}
\end{equation*}

\noindent
where $d_X(\,\,,\,\,)$ denotes the distance function on $X$ defined by the given hermitian metric.  We give natural orientationto $T_{\varepsilon}(c_{q-2})$.  Then,

\begin{equation}\label{e1.7}
\begin{split}
\lim_{\varepsilon\to 0}\int_{c_q-T_{\varepsilon}(c_{q-2})}\gamma(\alpha)&=
\lim_{\varepsilon\to 0}\int_{c_q-T_{\varepsilon}(c_{q-2})}d(\tilde{\alpha}\wedge\eta)\\
&= \lim_{\varepsilon\to 0}\int_{\partial T_{\varepsilon}(c_{q-2})}\tilde{\alpha}\wedge\eta \mbox{\hskip 20pt (by Stokes's Theorem)}
\end{split}
\end{equation}

\noindent
Using local coordinates $(z_1,\cdots,z_n,z_{n+1})$ on $X$ such that $Y$ 
is defined by $z_{n+1}=0$, $\tilde{\alpha}\wedge \eta$ and 
$\partial T_{\varepsilon}(c_{q-2})$ are locally written as 

\begin{eqnarray*}
\tilde{\alpha}\wedge \eta &=&\frac{1}{2\pi i}\tilde{\alpha}\wedge 
\frac{dz_{n+1}}{z_{n+1}}+(\mbox{{\it regular form\/}})\\
& & \\
\pm \partial T_{\varepsilon}(c_{q-2})&=&c_{q-2}\times 
\{~z_{n+1}\in \mathbb{C}~\vert ~\vert z_{n+1}\vert=\varepsilon~\}\\
& &\mbox{\hskip 70pt({\it with natural orientation)\/}}
\end{eqnarray*}

\noindent
Hence,
\begin{equation*}
\int_{\partial T_{\varepsilon}(c_{q-2})}\tilde{\alpha}\wedge \eta=
\pm\int_{c_{q-2}}
\tilde{\alpha}+\int_{\partial T_{\varepsilon}(c_{q-2})}(\mbox{{\it regular form)\/}},
\end{equation*}

\noindent
and since $\lim_{\varepsilon\to 0}\int_{\partial T_{\varepsilon}(c_{q-2})}(regular form)=0$, 

\begin{eqnarray}\label{e1.8}
\int_{\partial T_{\varepsilon}(c_{q-2})}\tilde{\alpha}\wedge \eta=
\pm\int_{c_{q-2}}\tilde{\alpha}
\end{eqnarray}

\noindent
From (\ref{e1.7}) and (\ref{e1.8}) it follows that the integral $\int_{c_q}\gamma(\alpha)$ converges and the equality in (\ref{e1.6}) holds as requied. 
\end{proof}
\smallskip

\begin{prop}\label{prop1.8}
We have the following commutative diagram:

\begin{equation}\label{e1.9}
\setlength{\unitlength}{0.7mm}
\begin{picture}(70,30)
\put(0,20){$H^q(X,\mathbb{C})$}
\put(22,20){\vector(3,-1){35}}
\put(35,18){$L$}
\put(10,18){\vector(0,-1){10}}
\put(5,13){$r^q$}
\put(0,2){$H^q(Y,\mathbb{C})$}
\put(23,3){\vector(1,0){25}}
\put(35,5){$G^q$}
\put(50,2){$H^{q+2}(X,\mathbb{C})$}
\put(60,0){\vector(0,-1){10}}
\put(63,-5){$r^{q+2}$}
\put(50,-16){$H^{q+2}(Y,\mathbb{C})$}
\put(15,0){\vector(3,-1){35}}
\put(30,-12){$L'$}
\end{picture}
\end{equation}

\bigskip
\bigskip
\bigskip
\bigskip
\noindent
where $L^\prime$ denotes the Hodge operator on $H^{\ast}(Y,\mathbb{C})$ 
associated to the fundamental form on $Y$, the restriction $\Omega_{\vert Y}$ of the fundamental form $\Omega$ to $Y$.
\end{prop}
\begin{proof}
We first show that the commutativity of the upper triangle. Let 
$\alpha$ be a closed $C^\infty$ $q$-form on $X$.  We denote by 
$[\alpha]\in H^{q}(X,\mathbb{C})$ its cohomology class.  Then,
\begin{eqnarray*}
(G^q\circ r^q)([\alpha])&=& [d(\alpha\wedge \eta)]\\
                         &=& [d\alpha\wedge \eta+(-1)^q\alpha\wedge d\eta]\\
                         &=& [d\eta\wedge \alpha].
\end{eqnarray*}
\noindent
 Now, we recall that $\omega:=d\eta$ is a closed (1.1)-form which represents 
 the first Chern class of the line bundle $[Y]$.  Hence, $\omega$ is 
 cohomologus to $\Omega$ in $H^{2}(X,\mathbb{C})$. From this it follows that
 \begin{equation*}
 [d\eta\wedge \alpha]=[\Omega\wedge \alpha]=L([\alpha]).
 \end{equation*}
 
 \noindent
 Thus $(G\circ r^\ast)([\alpha])=L([\alpha])$ as required.  Similarly, the 
 commutativity of the lower triangle can be proved.
 \end{proof}
 \smallskip
 
  We now return to the long exact sequence of cohomology (\ref{e1.3}).
By Theorem \ref{thm1.1}, Theorem \ref{thm1.2}, Theorem \ref{thm1.3}, Proposition \ref{prop1.8} and Grothendieck's theorem 
in \cite{Gro-1} which tells us (among other things) that 
$H^{q}(X-Y,\mathbb{C})=0$ for $q\geq n+2$, we can easily see that the long 
exact sequence of cohomology (\ref{e1.3}) breaks down into the short exact 
sequences as follows:

\begin{equation}\label{e1.10}
0\to H^q(X,\mathbb{C})\xrightarrow{r^q} H^q(X-Y,\mathbb{C})\to 0 
\qquad\mbox{for\,\, $0\leseq q\leseq 1$},
\end{equation}
\begin{equation}\label{e1.11}
\begin{split}
0\to & H^{q-2}(Y,\mathbb{C})\xrightarrow{G^{q-2}}H^{q}(X,\mathbb{C})
\xrightarrow{r^{q}}H^q(X-Y,\mathbb{C})\to 0 \\
&\mbox{\hskip 200pt for\,\, $2\leseq q\leseq n$},
\end{split}
\end{equation}
\begin{equation}\label{e1.12}
\begin{split}
&0 \to H^{n-1}(Y,\mathbb{C})\xrightarrow{G^{n-1}} H^{n+1}(X,\mathbb{C})
\xrightarrow{r^{n+1}} H^{n+1}(X-Y,\mathbb{C}) \\
 &\hskip 120pt \xrightarrow{R^{n+1}}H^{n}(Y,\mathbb{C})
\xrightarrow{G^n}H^{n+2}(X,\mathbb{C})\to 0,
\end{split}
\end{equation}
\begin{equation}\label{e1.13}
0\to H^q(Y,\mathbb{C})\xrightarrow{G^{q}}H^{q+2}(X,\mathbb{C})\to 0 \qquad 
n+1\leseq q\leseq 2n.
\end{equation}
 
  We now define the notions of {\it primitive cycles\/} and {\it finite cycles\/} on $X$ with respect to the prime section $Y$. 

\begin{defn}\label{def1.3}
A $q$-cycle $c_q$ on $X$ is defined to be {\it primitive\/} if its intersection cycle $c_q\cdot Y$ with $Y$ is zero in $H_{q-2}(Y,\mathbb{C})$.  
A $q$-cycle $c_q$ on $X$ is defined to be {\it finite\/} if its support is 
contained is contained in $X-Y$. 
\end{defn}
\smallskip

	We call homology classes of primitive cycles {\it primitive homology classes\/} and those of finite cycles {\it finite homology calsses\/}.   We denote the subspace of primitive (resp. finite) $q$-homology classes by $H_q(X,\mathbb{C})_0$ (resp. $H_q(X,\mathbb{C})_{f}$) and call it the {\it primitive $q$-homology group\/} of(resp. {\it finite $q$-homology groups\/} of $X$. Then by the definitions, 
\begin{eqnarray*}
H_q(X,\mathbb{C})_0&:=&{\rm Ker\/}\{~H_q(X,\mathbb{C})\xrightarrow{\cdot [Y]} 
H_{q-2}(Y,\mathbb{C})~\} \\
H_q(X,\mathbb{C})_f&:=&{\rm Im \/}\{~H_q(X-Y,\mathbb{C})
\xrightarrow{\iota_\ast} 
H_q(X,\mathbb{C})~\}.
\end{eqnarray*}
\begin{prop}\label{prop1.9}
Primitive $q$-cycles possibly exist on $X$ only for $q$ with $0\leseq q\leseq n+1$,\quad
\mbox{and}
\begin{equation*}
H_q(X,\mathbb{C})_0=H_q(X,\mathbb{C})_f \qquad \mbox{for\,\, $0\leseq q\leseq n+1$}.
\end{equation*}
\end{prop}
\begin{proof}
From the homology sequences dual to the cohomology sequences in (\ref{e1.10}) through (1.12) the assertion easily follows.
\end{proof}
\smallskip

	To state about the relation between primitive cohomology and 
homology groups, we introduce the notation for a subspace $S$ of $H^q(X,\mathbb{C})$ (resp. $H_q(X,\mathbb{C})$) as follows:
\begin{equation*}
Ann(S):=\{~[\alpha]\in H_q(X,\mathbb{C})~\vert ~\vert <[\omega],[\alpha]>=0 
\quad \mbox{for any $[\omega]\in S$}~\},
\end{equation*}

\noindent
where $<\,\,,\,\,>$ denotes the pairing between cohomology and homology. We call this the {\it annihilator\/} subspace of $H_q(X,\mathbb{C})$ by the 
subspace $S$.  
\smallskip

\begin{prop}\label{prop1.10}
\hskip 50pt
\begin{description}
\item[(i)] $H_q(X,\mathbb{C})\simeq H_q(X,\mathbb{C})_0$ \quad 
$(0\leseq q\leseq 1)$
\item[(ii)] $H_q(X,\mathbb{C})\simeq H_q(X,\mathbb{C})_0\oplus Ann(H_q(X,\mathbb{C}_0)$ \quad $(2\leseq q\leseq n+1)$
\end{description}
\end{prop}
\begin{proof} The assertion (i) follows from the definition of primitive homology.  We will now prove the assertion (ii).  By (i) of \ref{thm1.3} and Proposition 1.9, $G^{q-2}H^{q-2}(X,\mathbb{C})=LH^{q-2}(X,\mathbb{C})$.  Hence, 
by (ii) of Theorem 1.3, 
\begin{equation}\label{e1.14}
 H^{q-2}(X,\mathbb{C})\simeq G^{q-2}H^{q-2}(Y,\mathbb{C})\oplus  H^{q}(X,\mathbb{C})_0
\end{equation}

\noindent
Therefore, by duality
\begin{equation}\label{e1.15}
 H_{q}(X,\mathbb{C})\simeq Ann(G^{q-2}H^{q-2}(Y,\mathbb{C}))\oplus Ann(H^{q-2}(X,\mathbb{C})_0).
\end{equation}

\noindent
By considering the paring between the exact sequences of cohomology (\ref{e1.10}), (\ref{e1.11}) and their dual exact sequences of homology, 
\begin{equation}\label{e1.16}
Ann(G^{q-2}H^{q-2}(Y,\mathbb{C})\simeq \iota_\ast H_{q}(X-Y,\mathbb{C})
=H_{q}(X,\mathbb{C})_f
\end{equation}

\noindent
From (\ref{e1.15}), (\ref{e1.16}) and Proposition 1.11 follows the assertion (ii).
\end{proof}
\smallskip

\begin{prop}\label{prop1.11}
For\,\, $0\leseq q\leseq n+1$, 
$r^q: H^{q}(X,\mathbb{C})\to H^{q}(X-Y,\mathbb{C})$ is injective on the subspace $H^{q}(X,\mathbb{C})_0$ and 
\begin{equation*}
H^{q}(X,\mathbb{C})_0\simeq r^q H^{q}(X,\mathbb{C})\hookrightarrow H^{q}(X-Y,\mathbb{C}).
\end{equation*}
\end{prop}

\begin{proof} By the exactness of the cohomology sequences (\ref{e1.11}) and 
(\ref{e1.12}), $Im~G=ker~r^\ast$.  Hence the assertion follows from (\ref{e1.14}).\end{proof}
\smallskip

\begin{defn}\label{def1.4}
Cycles with compact support in $X-Y$ is defined to be {\it r{\'e}sidue cycle\/} 
if they bounds in $X$.  We call their homology classes {\it r{\'e}sidue homology classes\/}. 
\end{defn}

 We denote the subspace of $H^c_q(X-Y,\mathbb{C})$ comprising 
r{\'e}sidue homology classes by $H^c_q(X-Y,\mathbb{C})_{r\acute{e}s}$.  By the definition,
\begin{equation*}
H^c_q(X-Y,\mathbb{C})_{r\acute{e}s}=Ker~\{~H^c_q(X-Y,\mathbb{C})\xrightarrow{\iota_\ast}H_q(X,\mathbb{C})~\}.
\end{equation*}

\noindent
Actually, $H_q^c(X-Y,\mathbb{C})_{r\acute{e}s}\not=0$ only for\,\,$q=n+1$ and 
\begin{equation}\label{e1.17}
H^c_{n+1}(X-Y,\mathbb{C})_{r\acute{e}s}=\tau_{n}H_n(Y,\mathbb{C})
\end{equation}

\noindent
because of the exact homology sequence (\ref{e1.2}) which is dual to (\ref{e1.3}). 
\smallskip

\begin{prop}\label{prop1.12}
\begin{equation*}
r^{n+1} H^{n+1}(X,\mathbb{C})=Ann(H^c_{n+1}(X-Y,\mathbb{C})_{r\acute{e}s}).
\end{equation*}
\end{prop}

\begin{proof}
By considering the paring between the cohomology exact sequence (\ref{e1.12}) and its dual homology sequence, we have 
\begin{equation*}
r^{n+1}H^{n+1}(X,\mathbb{C})=Ann(\tau_{n} H_{n}(Y,\mathbb{C})).
\end{equation*}

\noindent
Hence the assertion follows from (\ref{e1.17}).
\end{proof}

We denote by $H^{q}(X,\mathbb{C})_0$ the primitive cohomology group with respect to the Hodge operator $L'$ on $Y$ which is associated to $\Omega_{\vert Y}$, the restriction of the fundamental form $\Omega$ on $X$.  We are now going to 
discuss the primitive cohomology and homology of $Y$.  For use later we wish to make clear the relation between the image of the map $R^{n+1}:H^{n+1}(X-Y,\mathbb{C})\to H^{n}(Y,\mathbb{C}))$ in the exact sequence (\ref{e1.12}) and the primitive cohomology group $H^{n}(Y,\mathbb{C})_0$.  The result is as follows:

\begin{lem}\label{lem1.13}
The restriction map $r^n:H^{n}(X,\mathbb{C})\to H^{n}(Y,\mathbb{C})$, which is injective by the Weak Lefshetz Thorem, give rise to an isomorphism from 
$H^{n}(Y,\mathbb{C})_0$ into $H^{n}(Y,\mathbb{C})_0$ and 
\begin{equation*}
r^n(H^{n}(X,\mathbb{C}))\cap H^{n}(Y,\mathbb{C})_0=r^n(H^{n}(X,\mathbb{C})_0)
\end{equation*}
\end{lem}

\begin{proof} By the definition of primitive cohomology, the isomorphism in 
(\ref{e1.13}) for $n+2$, and \ref{thm1.3}, (ii), we have the following commutative diagram of exact sequences:
\begin{equation}\label{e1.18}
\begin{CD}
& & & & 0 & & \\
& & & & @VVV & & \\
0 @>>> H^{n}(X,\mathbb{C})_0 @>>> H^{n}(X,\mathbb{C}) @>\text{$L^2$}>> H^{n+4}(X,\mathbb{C})\\
& & & & @V\text{$r^n$}VV  @A\text{$\simeq$}A\text{$G^{n+2}$}A \\
0 @>>> H^{n}(Y,\mathbb{C})_0 @>>> H^{n}(Y,\mathbb{C}) @>\text{$L^{\prime}$}>> H^{n+2}(Y,\mathbb{C}) \\
\end{CD}
\end{equation}
\bigskip

\noindent
From this we infer that $r^n(H^{n}(Y,\mathbb{C}))_0\hookrightarrow H^{n}(Y,\mathbb{C})_0$ and $r_{\vert H^{n}(Y,\mathbb{C})_0}^{n}(H^{n}(Y,\mathbb{C})_0\to H^{n}(Y,\mathbb{C})$ is an isomorphism into.  To show the latter part, we 
consider the following Lefshetz decompositions of $H^{n}(Y,\mathbb{C})$ and 
$H^{n}(Y,\mathbb{C})$:

\begin{equation}\label{e1.19}
H^{n}(X,\mathbb{C})=H^{n}(X,\mathbb{C})_0\oplus LH^{n-2}(X,\mathbb{C}),
\end{equation}
\begin{equation}\label{e1.20}
H^{n}(Y,\mathbb{C})=H^{n}(Y,\mathbb{C})_0\oplus L^{\prime}H^{n-2}(Y,\mathbb{C}).\end{equation}
\smallskip

\noindent
Note that, since $r^{n-2}:H^{n-2}(X,\mathbb{C})\to H^{n-2}(Y,\mathbb{C})$ is an isomorphism by the Weak Lefshetz Theorem, $r^n:H^{n}(X,\mathbb{C})\to H^{n}(Y,\mathbb{C})$ maps $LH^{n-2}(X,\mathbb{C})$ onto $L^{\prime}H^{n-2}(X,\mathbb{C})$ isomorphically.  The inclusion

\begin{equation}\label{e1.21}
 r^nH^{n}(X,\mathbb{C})_0 \hookrightarrow r^nH^{n}(X,\mathbb{C})\cap H^{n}(Y,\mathbb{C})_0
\end{equation}
\smallskip

\noindent
is obvious, since $r^nH^{n}(X,\mathbb{C})_0\hookrightarrow H^{n}(Y,\mathbb{C})_0$ as has been proved just above.  We will prove the reverse inclusion.  
Given $x\in r^nH^{n}(X,\mathbb{C})\cap H^{n}(Y,\mathbb{C})_0$, there 
exists a $y\in H^{n}(X,\mathbb{C})$ with $r^n(y)=x$.  We write $y$ as 
$y=y_1+y_2$ where $y_1\in H^{n}(X,\mathbb{C})_0$ and $y_2\in LH^{n-2}(X,\mathbb{C})$.  Then $x=r^n(y)=r^n(y_1)+r^n(y_2)$, and $r^n(y_1)\in 
H^{n}(Y,\mathbb{C})_0$, $r_n^\ast(y_2)\in L^{\prime}H^{n-2}(Y,\mathbb{C})_0$.  
Hence
 
\begin{equation*}
x-r^n(y_1)=r^n(y_2)\in H^{n}(Y,\mathbb{C})_0\cap L^{\prime}H^{n-2}(Y,\mathbb{C})_0=0.
\end{equation*}
\smallskip

\noindent
Thus $r^n(y_2)=0$, from which $y_2=0$ follows since $r^n$ maps 
$LH^{n-2}(X,\mathbb{C})$ onto $ L^{\prime}H^{n-2}(Y,\mathbb{C})$ isomorphically. Hence $x=r^n(y_1)$.  This shows that 
\begin{equation}\label{e1.22}
r^nH^{n}(X,\mathbb{C})\cap H^{n}(Y,\mathbb{C})_0\hookrightarrow r^nH^{n}(X,\mathbb{C})_0
\end{equation}

\noindent
By (\ref{e1.21}) and (\ref{e1.22}), $r^nH^{n}(X,\mathbb{C})\cap H^{n}(Y,\mathbb{C})_0=r^nH^{n}(X,\mathbb{C})_0$ as required. 
\end{proof}
\smallskip 

\begin{lem}\label{lem1.14}
There is an exact sequence 
\begin{equation}\label{e1.23}
0\to LH^{n}(X,\mathbb{C})_0\to H^{n+2}(X,\mathbb{C})\xrightarrow{r^{n+2}}
H^{n+2}(Y,\mathbb{C})\to 0.
\end{equation}
\end{lem}

\begin{proof}
To see the surjectivity of $r_{n+2}^\ast$, we consider the following commutative diagram:

\begin{equation*}
\begin{CD}
H^{n+2}(X,\mathbb{C}) @>\text{$r^{n+2}$}>> H^{n+2}(Y,\mathbb{C}) & &\\
 @A\text{$L$}A\text{$\simeq$}A   @A\text{$\simeq$}A\text{$L^{\prime 2}$}A & &\\
H^{n}(X,\mathbb{C}) @<<\text{$G^{n-2}$}< H^{n-2}(Y,\mathbb{C}) @<<<0
\end{CD}
\end{equation*} 
\smallskip

\noindent
the commutativity of which follows from the description of the Gysin map $G^{n-2}$ using differential forms (\ref{prop1.8}). From this diagram the surjectivity of $r^{n+2}$ follows, since $L^{\prime 2}$ is an isomorphism (Hard Lefshetz for $Y$). The injectivity of $L: H^{n}(X,\mathbb{C})_0\to H^{n+2}(X,\mathbb{C})$ follows from the fact that $L: H^{n}(X,\mathbb{C})_0\to H^{n+2}(X,\mathbb{C})$ is an isomorphism (Hard Lefshetz for $X$). \par
  To prove the exactness at the term $H^{n+2}(X,\mathbb{C})$, we consider the following commutative diagram: 
\begin{equation}\label{e1.24}
\begin{CD}
& & 0 & & 0 & & & & \\
& & @VVV  @VVV & & & & \\
& & H^{n}(X,\mathbb{C})_0 @>\text{{\it inclusion}}>> H^{n}(X,\mathbb{C})_0 @>\text{$\simeq$}>\text{$L$}>H^{n+2}(X,\mathbb{C}) & & \\
& & @V\text{\scriptsize{(Lemma 1.13)}}V\text{$r^n$}V   @VV\text{$r^n$}V   @VV\text{$r^{n+2}$}V & & \\
0 @>>> H^{n}(Y,\mathbb{C})_0 @>>> H^{n}(Y,\mathbb{C}) @>\text{$L^\prime$}>> 
H^{n+2}(Y,\mathbb{C}) @>>> 0 \text{\quad (exact)} \\
& & & & & & @VVV & & \\
& & & & & & 0 & &
\end{CD}
\end{equation}  
\smallskip

\noindent
By this diagram we can easily see $LH^{n}(X,\mathbb{C})_0\subset Ker~r^n$.  We will prove the converse inclusion by casing the diagram (\ref{e1.24}).  Given $x\in Ker~r^n$, there exists a $y\in H^{n}(X,\mathbb{C})$ with $L(y)=x$.  Then $L^\prime(r^n(y))=r^{n+2}(L(y))=r^{n+2}(x)=0$, hence 
$r^n(y)\in r^nH^{n}(X,\mathbb{C})\cap H^{n}(X,\mathbb{C})_0$.  We should now recall that 
$r^nH^{n}(X,\mathbb{C})\cap H^{n}(Y,\mathbb{C})_0=r^nH^{n}(X,\mathbb{C})_0$ (Lemma 1.13) and $r^n$ is injective.  This implies $y\in H^{n}(X,\mathbb{C})_0$, that is, $x=L(y)\in LH^{n}(X,\mathbb{C})_0$, which means 
$Ker~r^n\subset LH^{n}(X,\mathbb{C})_0$.  Consequently, we conclude 
$Ker~r^n=LH^{n}(X,\mathbb{C})_0$ as requied.  
\end{proof}
\smallskip

\begin{thm}\label{thm1.15}
\begin{equation*}
H^{n}(Y,\mathbb{C})_0=R^{n+1}(H^{n+1}(X-Y,\mathbb{C})\oplus r^n(H^{n}(X,\mathbb{C})_0)
\end{equation*}
\end{thm}

\begin{proof}
Let us consider the Lefshetz decompositions of $H^{n}(Y,\mathbb{C})_0$ and 
$H^{n+2}(X,\mathbb{C})$:

\begin{eqnarray*}
& &H^{n}(Y,\mathbb{C})=H^{n}(Y,\mathbb{C})_0\oplus L^{\prime}H^{n-2}(Y,\mathbb{C})\\
& &H^{n+2}(X,\mathbb{C})=LH^{n}(X,\mathbb{C})_0\oplus L^{2}H^{n-2}(X,\mathbb{C}). 
\end{eqnarray*}

{\bf Claim:} Concerning the Gysin map $G^n:H^{n}(Y,\mathbb{C})\to H^{n+2}(X,\mathbb{C})$, we have 

\begin{description}
\item[(a)] $G^n(L^{\prime}H^{n-2}(Y,\mathbb{C})\subset L^2H^{n-2}(X,\mathbb{C})$\quad $G^n$ maps $L^{\prime}H^{n-2}(Y,\mathbb{C})$ onto \\
$L^{2}H^{n-2}(X,\mathbb{C})$ isomorphically, 
\item[(b)] $G^n(H^{n}(Y,\mathbb{C})_0)=LH^{n}(X,\mathbb{C})_0$, \quad and 
\item[(c)] $Ker~G^n\subset H^{n}(Y,\mathbb{C})_0$. 
\end{description}

{\it Proof of {\rm (a)}}:  By Proposition \ref{prop1.9}, we have the following commutative diagram: 

\begin{equation*}
\setlength{\unitlength}{0.7mm}
\begin{picture}(130,30)(0,-5)
\put(0,22){$H^{n-2}(X,\mathbb{C})$}
\put(23,18){\vector(3,-1){35}}
\put(35,16){$L$}
\put(10,19){\vector(0,-1){11}}
\put(-2,13){$r^{n-2}$}
\put(13,13){$\simeq$}
\put(0,2){$H^{n-2}(Y,\mathbb{C})$}
\put(25,3){\vector(1,0){25}}
\put(33,5){$G^{n-2}$}
\put(50,2){$H^{n}(X,\mathbb{C})$}
\put(60,0){\vector(0,-1){10}}
\put(62,-6){$r^{n}$}
\put(50,-16){$H^{n}(Y,\mathbb{C})$}
\put(15,0){\vector(3,-1){35}}
\put(30,-12){$L'$}
\put(74,1){\vector(3,-1){35}}
\put(91,-3){$L$}
\put(79,-13){$G^n$}
\put(73,-15){\vector(1,0){25}}
\put(100,-16){$H^{n+2}(X,\mathbb{C})$}
\end{picture}
\end{equation*}
\bigskip
\bigskip

\noindent
where $r^{n-2}:H^{n-2}(X,\mathbb{C})\to H^{n-2}(Y,\mathbb{C})$ is 
an isomorphism by the Weak Lefshetz Theorem.  From this diagram 
$G^n(L^{\prime}H^{n-2}(Y,\mathbb{C})\subset L^2H^{n-2}(X,\mathbb{C})$ follows. 
The fact that $G^n$ maps $L^{\prime}H^{n-2}(Y,\mathbb{C})$ onto $L^{2}H^{n-2}(X,\mathbb{C})$ isomorphically is proved as follows: Since $L:H^{n-2}(X,\mathbb{C})
\to H^{n}(X,\mathbb{C})$ is injective, and since $L:H^{n}(X,\mathbb{C})
\to H^{n+2}(X,\mathbb{C})$ is an isomorphism (Hard Lefshetz Theorem), 
$L^2:H^{n-2}(X,\mathbb{C})\to L^2H^{n+2}(X,\mathbb{C})$ is an isomorphism. 
Besides, since $L^{\prime}:H^{n-2}(Y,\mathbb{C})\to H^{n}(Y,\mathbb{C})$ is 
injective, $L^{\prime}:H^{n-2}(Y,\mathbb{C})\to L^{\prime}H^{n-2}(Y,\mathbb{C})$ is also an isomorphism.  Therefore, taking into account that $r^{n-2}:H^{n-2}(X,\mathbb{C})\to H^{n-2}(Y,\mathbb{C})$ is an isomorphism, we conclude that the Gysin map $G^n$ maps  $L^{\prime}H^{n-2}(Y,\mathbb{C})$ onto $L^{2}H^{n-2}(X,\mathbb{C})$ isomorphically. \par
\smallskip

 {\it Proof of {\rm (b)}}: Combining (\ref{e1.9}) for $q=n$, Proposition \ref{prop1.9},
 (\ref{e1.23}) and (\ref{lem1.13}), we have the following commutative diagram:
\begin{equation*}
\setlength{\unitlength}{0.7mm}
\begin{picture}(150,60)(0,-30)
\put(0,0){$0$}
\put(4,2){\vector(1,0){10}}
\put(15,0){$H^{n}(Y,\mathbb{C})_0$}
\put(40,2){\vector(1,0){10}}
\put(51,0){$H^{n}(Y,\mathbb{C})$}
\put(75,2){\vector(1,0){11}}
\put(79,4){$L^\prime$}
\put(87,0){$H^{n+2}(Y,\mathbb{C})$}
\put(114,2){\vector(1,0){10}}
\put(125,0){$0\quad\text{(exact)}$}
\put(98,6){\vector(3,2){16}}
\put(116,17){$0$}
\put(119,23){(exact)}
\put(61,-2){\vector(0,-1){10}}
\put(64,-8){$G^n$}
\put(51,-18){$H^{n+2}(X,\mathbb{C})$}
\put(73,-12){\vector(3,2){16}}
\put(61,-19){\vector(0,-1){10}}
\put(60,-34){$0$}
\put(15,-34){$LH^{n}(X,\mathbb{C})_0$}
\put(42,-30){\vector(3,2){16}}
\put(11,-49){$0$}
\put(15,-46){\vector(3,2){16}}
\end{picture}
\end{equation*}
\bigskip
\bigskip
\bigskip

\noindent
From this it follows 
\begin{equation*}
G^n(H^{n}(Y,\mathbb{C})_0)\subset Ker~r^{\ast}_{n+2}=LH^{n}(X,\mathbb{C})_0
\end{equation*}

\noindent
Actually, they coincides with each other, since $G^n$ is surjective and (a) holds.\smallskip

{\it Proof of {\rm (c)}:} Let $x\in Ker~G^n$. We write it as $x=x_1+x_2$, where $x_1\in H^{n}(Y,\mathbb{C})_0$ and $x_2\in L^{\prime}H^{n-2}(Y,\mathbb{C})_0$.  Then $G^n(x)=G^n(x_1)+G^n(x_2)=0$, and by (a) and (b), $G^n(x_1)\in LH^{n}(X,\mathbb{C})_0$ and $G^n(x_2)\in L^2H^{n-2}(X,\mathbb{C})_0$.  Hence $G^n(x_2)=-G^n(x_1)\in LH^{n}(X,\mathbb{C})_0\cap L^2H^{n-2}(X,\mathbb{C})_0=0$.  Thus $G^n(x_2)=0$, whence $x_2=0$. This is because $G^n$ maps $L^{\prime}H^{n-2}(Y,\mathbb{C})$ onto $L^{2}H^{n-2}(Y,\mathbb{C})$ isomorphically.  Therefore 
$x=x_1\in H^{n}(Y,\mathbb{C})_0$, which means $Ker~G^n\subset H^{n}(Y,\mathbb{C})_0$. \par
\hfill{q.e.d. for the Claim.}
\bigskip

 Now we can easily deduce the Proposition.  In fact, by Lemma \ref{lem1.13} and the claim (a), (b) (c) above, we have the following commutative diagram:\par
\begin{equation*}
\setlength{\unitlength}{0.7mm}
\begin{picture}(130,70)(0,-32)
\put(60,34){$0$}
\put(62,32){\vector(0,-1){10}}
\put(55,18){$Ker~G^n$}
\put(62,17){\vector(0,-1){10}}
\put(0,2){$0$}
\put(5,4){\vector(1,0){10}}
\put(16,2){$H^{n}(X,\mathbb{C})_0$}
\put(45,5){$r^n$}
\put(41,4){\vector(1,0){13}}
\put(55,2){$H^{n}(Y,\mathbb{C})_0$}
\put(40,0){\vector(2,-1){18}}
\put(35,-6){$L$}
\put(49,-3){$\simeq$}
\put(62,0){\vector(0,-1){10}}
\put(64,-5){$G^n$}
\put(56,-15){$LH^{n}(Y,\mathbb{C})_0$}
\put(62,-17){\vector(0,-1){10}}
\put(61,-32){$0$,}
\end{picture}
\end{equation*}

\noindent
which implies

\begin{equation*}
H^{n}(Y,\mathbb{C})_0\simeq Ker~G^n\oplus r^{n}(H^{n}(X,\mathbb{C})_0).
\end{equation*}

\noindent
Here, recall that $Ker~G^n=Im~R^{n+1}$ by (1.12), then we are done.
\end{proof}
\smallskip

 We wish to identify the subspace of $H^{n}(Y,\mathbb{C})_0$ which is dual to 
 $Im~R^{n+1}$.  For this puopose we need to introduce the following notion.
 
\begin{defn}\label{defn1.5}
Cycles in $Y$ is defined to be {\it vanishing cycles\/} with respect to $X$ if 
they bound in $X$.  We call their homology classes {\it vanishing homology 
classes\/}.
\end{defn}
\smallskip

 We denote the subspace $H_{q}(Y,\mathbb{C})$ comprising vanishing homology 
 classes by $H_{q}(Y,\mathbb{C})_{v}$.  Note that $H_{q}(Y,\mathbb{C})_v$ may 
 not be zero only if $q=n$.  
 \smallskip
 
 \begin{prop}\label{prop1.16}
 $H_{q}(Y,\mathbb{C})_v$ is included in $H_{n}(Y,\mathbb{C})_0$ and 
 \begin{equation*}
 H_{n}(Y,\mathbb{C})=H_{n}(Y,\mathbb{C})_v\oplus Ann(Im~R^{n+1})
 \end{equation*}
 
\noindent
or, equivalently 
\begin{equation*}
 H_{n}(Y,\mathbb{C})_0=H_{n}(Y,\mathbb{C})_v\oplus [Ann(Im~R^{n+1})\cap H_{n}(Y,\mathbb{C})_0]
\end{equation*}
\end{prop}

\begin{proof}

By virture of Theorem \ref{thm1.15}, it suffices to show that 
\begin{equation*}
 H_{n}(Y,\mathbb{C})_0\cap Ann(r^{n}H^{n}(X,\mathbb{C})_0)= H_{n}(Y,\mathbb{C})_v.
\end{equation*} 

\noindent
The inclusion $H_{n}(Y,\mathbb{C})_v\subset Ann(r^{n}H^{n}(X,\mathbb{C})_0)$ is trivial.  To see that $H_{n}(Y,\mathbb{C})_v\subset H_{n}(Y,\mathbb{C})_0$, 
consider the following diagram: 
\begin{equation}\label{e1.25}
\begin{CD}
H_{n}(Y,\mathbb{C}) @>\text{$\cdot[Z]$}>>  H_{n-2}(Z,\mathbb{C})_v \\
  @V\text{$\iota_n$}VV  @V\text{$\simeq$}VV\text{$\iota_{n-2}$}\\
H_{n}(X,\mathbb{C}) @>>\text{$\cdot[Y]$}>  H_{n-2}(Y,\mathbb{C})_v,
\end{CD}
\end{equation}

\noindent
where $Z$ is the intersection of a generic member $\vert Y\vert$ (linear sysytem of effective divisors which are linearly equivalent to $Y$) with $Y$, which is a non-singular, irreducible hypersurface of $Y$ and for which $c_1([Z])\sim \Omega_{\vert Y}$ (cohomologous), where $\iota^n$ (resp. $\iota^{n-2}$) is the homomorphism induced by the inclusion map $\iota:Y\hookrightarrow X$ 
(resp. $\iota:Z\hookrightarrow Y$, and where $\cdot [Z]$ (resp. $[Y]$) is the 
map which assignes each $n$-cycle in $Y$ (resp. $X$) to its intersection cycle with $Z$ (resp. $Y$). By the diagram (\ref{e1.25}), 
$H_{n}(Y,\mathbb{C})_v\hookrightarrow Ker~(\cdot [Z])$.  Meanwhile, 
$Ker~(\cdot [Z])=H_{n}(Y,\mathbb{C})_0$ by definition.  Thus we have 
$H_{n}(Y,\mathbb{C})_v\hookrightarrow H_{n}(Y,\mathbb{C})_0$.  Hence 
\begin{equation}\label{e1.26}
H_{n}(Y,\mathbb{C})_v\hookrightarrow H_{n}(Y,\mathbb{C})_0\cap Ann(r^{n}H^{n}(X,\mathbb{C})_0).
\end{equation}

\noindent
Next we will prove the converse inclusion.  It suffices to show that if 
$[\gamma]\in H_{n}(Y,\mathbb{C})_0\cap Ann(r^{n}H^{n}(X,\mathbb{C})_0)$, then $\int_{\gamma}\omega=0$ for any $[\omega]\in H^{n}(X,\mathbb{C})$.  To see this, we use the Lefshetz decomposition
\begin{equation*}
H^{n}(X,\mathbb{C})=H^{n}(Y,\mathbb{C})_0\oplus LH^{n-2}(X,\mathbb{C}).
\end{equation*}

\noindent
Assume $[\gamma]\in H_{n}(Y,\mathbb{C})_0\cap Ann(r^{n}H^{n}(X,\mathbb{C})_0)$.  Then $\int_{\gamma}\omega=0$ for any $[\omega]\in  H^{n}(X,\mathbb{C})_0$, 
and for any $[\Omega\wedge \omega^\prime]\in LH^{n-2}(X,\mathbb{C})$ ($\omega^\prime]\in H^{n-2}(X,\mathbb{C})$,
\begin{equation*}
\int_{\gamma}\Omega\wedge \omega^\prime=\int_{[\gamma\cdot Y]}\omega^\prime=0
\end{equation*}

\noindent
since $[\gamma\cdot Y]=0$ by the assumption.  Thus $\int_{\gamma}\omega=0$ for 
any $[\omega]\in H^{n}(X,\mathbb{C})$ if $[\gamma]\in H_{n}(Y,\mathbb{C})_0\cap Ann(r^{n}H^{n}(X,\mathbb{C})_0)$.  This implies 
\begin{equation}\label{e1.27}
H_{n}(Y,\mathbb{C})_v\hookleftarrow H_{n}(Y,\mathbb{C})_0\cap Ann(r^{n}H^{n}(X,\mathbb{C})_0)
\end{equation}

\noindent
By (\ref{e1.26}) and (\ref{e1.27}), $H_{n}(Y,\mathbb{C})_v= H_{n}(Y,\mathbb{C})_0\cap Ann(r^{n}H^{n}(X,\mathbb{C})_0)$ as requied.
\end{proof}
\bigskip
\section{Rational De Rham groups  
of an algebraic manifold and Integrals of the second kind on it}

As in \S 1 we let $X$ be a non-singular irreducible algebraic variety of 
dimension $n+1$ embedded in a higher dimensional complex projecyive space $\mathbb{P}^N(\mathbb{C})$ and $Y$ a generic hyperplane section of $X$.  By a {\it meromorphic $q$-form on $X$\/} we shall mean an exterior differential form $\omega$ of degree $q$, which has the form 
\begin{equation*}
\omega=\sum f_{i_1i_2\cdots i_q}dz_{i_1}\wedge dz_{i_2}\wedge\cdots\wedge dz_{i_q}
\end{equation*}

\noindent
where $(z_1,\cdots,z_{n+1})$ is a complex analytic local coordinate 
system on $X$ and $f_{i_1i_2\cdots i_q}{'s}$ are meromorphic functions of the variables $(z_1,\cdots,z_{n+1})$.  We denote by $\Omega^q_X(kY)$ the sheaf of 
germs of meromorphic $q$-forms having poles of order $k$ (at most) along $Y$ 
as their only sngularities.  The direct limit of the sheaves $\Omega^q_X(kY)$ a $k\to \infty$ we denote by $\Omega^q_X(\ast Y)$.  It is just the sheaf of germs of meromorphic $q$-forms with poles of arbitrary order along $Y$.  We put $\Omega^{\cdot}_X(\ast Y):=\sum \Omega^q_X(\ast Y)$, which forms a complex of sheaves with respect to the exterior derivative $d$. We define 
\begin{equation*}
\Phi^q(kY):=Ker~\{~\Omega^q_X(kY)\xrightarrow{d}\Omega^{q+1}_X((k+1)Y)~\}
\end{equation*}

\noindent
and call it the sheaf of germs of {\it closed meromorphic $q$-forms\/} having 
 poles of order $k$ (at most) along $Y$ as their only singularities.  We define the sheaf $\Omega^q_X(\log Y)$ to be the subsheaf of $\Omega^q_X(\ast Y)$ consisting of 
the germs of such local meromorphic $q$-forms that both of $f\omega$ and 
$df\wedge \omega$ are holomorphic if $f$ is a local holomorphic defining 
equation of $Y$.  If $g=0$ is another defining equation of $Y$, then $g=uf$ where $u$ is a non-vanishing local holomorphic function and the relation 
$g\omega=uf\omega$, $dg\wedge \omega=udf\wedge \omega+fdu\wedge \omega$ shows that $\Omega^q_X(\log Y)$ is well-defined.  We call the sheaf of germs of meromorphic $q$-forms having {\it logarithmic poles (at most) along $Y$\/} as thier only singularities. The reason for this naming is that a meromorphic $q$-form $\omega$ ($q\geq 1$) has logarithmic poles (at most) along $Y$ as its only singularities if and only if $\omega$ is locally written as 
\begin{equation*}
\omega=\varphi\wedge \frac{df}{f}+\psi,
\end{equation*}

\noindent
where $\phi$, $\psi$ are holomorphic forms and $f=0$ is a local holomorphic equation of $Y$.  The following lemma is fundametal for calculations in the subsequel.  

\begin{lem}\label{lem2.1}
\hskip 5pt \par
\vskip 3pt
\noindent
\hskip 3pt {\rm (i)} The following sheaf sequences are exact:\par
\vskip 5pt
\noindent
\hskip 5pt {\rm (a)} \quad $0\to \Phi^{q-1}((k-1)Y)\to \Omega^{q-1}_X((k-1)Y)\xrightarrow{d}
\Phi^{q}(kY)\to 0 \quad (q\geq 2, k\geq 2)$ \par
\vskip 5pt
\noindent
\hskip 5pt {\rm (b)}\qquad $0\to \Phi^{q-1}(Y)\to \Omega^{q-1}_X(\log Y)\xrightarrow{d}
\Phi^{q}(Y)\to 0 \quad (q\geq 2)$ \par
\vskip 10pt
\noindent
\hskip 3pt {\rm (ii)} There exist naturally the following exact sequences of sheaves: \par
\vskip 5pt
\noindent
\hskip 5pt{\rm (c)}\hskip 15pt $0\to \mathbb{C}_X \to {\cal O}((k-1)Y)\xrightarrow{d}\Phi^{1}(kY)\xrightarrow{\alpha} \mathbb{C}_Y \to 0 \quad (k\geq 1)$\par
\vskip 5pt
\noindent
\hskip 5pt {\rm (d)}\hskip 20pt $0\to \Omega^{q}_X(Y)\to \Omega^{q}_X(\log Y)\xrightarrow{R}
\Omega^{q-1}_X(Y)\to 0 \quad (q\geq 1)$ \par
\vskip 5pt
\noindent
\hskip 5pt {\rm (e)}\hskip 30pt $0\to \Phi^{q}_X\to \Phi^{q}_X(Y)\xrightarrow{R}
\Phi^{q-1}_Y\to 0 \quad (q\geq 1)$
\end{lem}

\begin{proof} We take a local coordinate system $(z_1,\cdots,z_n,w)$ on $X$ such that $Y$ is defined by $w=0$.  First, we prove for all pairs of integers $(q,k)$ with $q\geq 1,\,\,k\geq 1$ that if $\varphi$ is a local holomorphic section of the sfeaf $\Phi^q(kY)$, then $\varphi$ is written as 
\begin{equation}\label{e2.1}
\varphi=\frac{A\wedge dw}{w^k}+\frac{B}{w^{k-1}}
\end{equation}

\noindent
where $A$, $B$ are holomorphic and involve only $dz_1,\cdots,dz_{n}$.  In fact, as to such $\varphi$, since $w^k\varphi$ is holomorphic, we may write 
\begin{equation*}
\varphi=\frac{A\wedge dw}{w^k}+\frac{B^\prime}{w^{k}}
\end{equation*}

\noindent
where $A$, $B^\prime$ are holomorphic and do not involve $dw$.  Since $\varphi$ is closed, 
\begin{equation*}
d\varphi=\frac{dA\wedge dw}{w^k}+\frac{dB^\prime}{w^{k}}+(-k)\frac{dw\wedge B^\prime}{w^{k+1}}=0
\end{equation*}

\noindent
so that $B:=B^\prime/w$ is holomorphic.  Hence we have locally the expression in (\ref{e2.1}) as required.  Now we prove the exactness of (i)-(a) and (i)-(b).  For a local holomorphic section $\varphi$ of $\Phi^q(kY)$ (~$q\geq 1,~k\geq 2$ and $q\geq 2,~k=1$), we take such an expression as in (\ref{e2.1}).  If $k\geq 2$, letting $\psi_1=-(1/(k-1))(A/w^{k-1})$, $\varphi-d\psi_1$ is a local section of $\Phi^q(k-1)$.  Repeating this argument, we may find a local section $\psi$ of $\Omega^{q-1}_X((k-1)Y)$ such that $\varphi-d\psi$ is a section of $\Phi_X^q(Y)$.  Thus 
\begin{equation*}
\varphi-d\psi=E\wedge \frac{dw}{w}+F,
\end{equation*}

\noindent
where $E$, $F$ are holomorphic and involve only $dz_1\cdots,dz_{n}$. We express $E$ as follows:
\begin{equation*}
E=E_0(z)+wE_1(z,w)
\end{equation*}

\noindent
where $E_0(z)$ does not involve $w$.  Then,
\begin{equation*}
\varphi-d\psi=E_0(z)\wedge \frac{dw}{w}+F_0
\end{equation*}

\noindent
where $F_0=E_1+F$. Since $d(Edw/w+F)=0$, $d_zE_0(z)dw/w+dF_0=0$.  Hence 
$d_zE_0(z)dw+wdF_0=0$.  From this it follows that $d_zE_0(z)=0$, $dF_0=0$. Therefore, there exist $D(z)$ and $G$ such that $d_zD=E_0$ and $dG=F_0$, and so 
\begin{equation*}
d(D\frac{dw}{w}+G)=E_0\wedge \frac{dw}{w}+F_0.
\end{equation*}

\noindent
Hence,
\begin{equation}\label{e2.2}
\varphi=d(\psi+D\wedge \frac{dw}{w}+G),
\end{equation}

\noindent
namely, $\varphi$ is a derived form. This shows the exactness of the sequence (i)-(a). If $k=1$, then $\psi$ does not appear in the expression of $\varphi$ in (\ref{e2.2}).  This shows the exactness of (i)-(b).\par
\smallskip

  Next we prove the exactness of the sequence (ii)-(c).  If $\varphi$ is a local section $\Phi^1(Y)$, then it is written as 
\begin{equation*}
\varphi=A\wedge \frac{dw}{w}+B,
\end{equation*}

\noindent
where $A$ is a holomorphic function and $B$ is a holomorphic $1$-form, involving only $dz_1\cdots,dz_{n}$ (cf. (\ref{e2.1})).  Writting $A$ as 
\begin{equation*}
A(z,w)=A_0(z)+wA_1(z,w),
\end{equation*}

\noindent
where $A_0(z)$ is a function of $z_1,\cdots,z_{n}$, we have
\begin{equation*}
\varphi=A_0(z)\wedge \frac{dw}{w}+B_0,
\end{equation*}

\noindent
where $B_0=A_1(z,w)dw+B$.  Since 
\begin{equation*}
d\varphi=\frac{d_zA_0(z)\wedge dw}{w}+dB_0=0,
\end{equation*}

\noindent
we have
\begin{equation*}
d_zA_0(z)dw+w~dB_0=0,
\end{equation*}

\noindent
Hence $d_zA_0(z)=dB_0=0$.  From these it follows that $A_0(z)$ is constant and 
$B_0=dC$ for some holomorphic function $C(z,w)$.  Thus $\varphi$ is written as
\begin{equation*}
\varphi=A_0\wedge \frac{dw}{w}+dC,
\end{equation*}

\noindent
This means $\Phi^1_X(Y)/d\Omega^0_X$ is locally a constant sheaf.  At each point $y\in Y$, we take $[(1/2\pi i)dw/w]_y$, the class of $(\Phi^1_X(Y)/d\Omega^0_X)_y$ 
determined by $(1/2\pi i)dw/w$, as a generator of $(\Phi^1_X(Y)/d\Omega^0_X)_y$.  We can easily see that the class $[(1/2\pi i)dw/w]_y$ is uniquely determined,
not depending on the choice of a local defining equation of $Y$.  We denote by 
$\alpha_y:(\Phi^1_X(Y)/d\Omega^0_X)_y\to \mathbb{C}_{Y,y}$ defined by 
\begin{equation*}
 \Big[%
 \frac{1}{2\pi i}\frac{dw}{w}
 \Big]%
 _y \to 1_{Y,y}
\end{equation*}

\noindent
at each point $y\in Y$, which gives rise to a well-defined sheaf homomorphism 
$\alpha:\Phi^1_X(Y)/d\Omega^0_X\to \mathbb{C}_{Y}$ as easily seen. The surjectivity of the map $\alpha$ and that the kernel of the homomorphism 
$d:\Omega^0_X\to \Phi^1_X(Y)$ coincides with $\mathbb{C}_X$ is obvious.  The sheaf 
homomorphism $R^q:\Omega_X^q(\log Y)\to \Omega_Y^{q-1}$, which we call 
{\it R\'esidues map\/} is defined as follows (resp. $R:\Phi_X^q(Y)\to \Phi_Y^{q-1}$): A local cross-section $\varphi$ of the sheaf $\Omega_X^q(\log Y)$ (resp. $\Phi_X^q(Y)$) is written as 
\begin{equation*}
\omega=\varphi\wedge \frac{dw}{w}+\psi,
\end{equation*}

\noindent
where $\varphi$ is a holomorphic $(q-1)$-form and $\psi$ is a holomorphic $q$-form, involving only $dz_1,\cdots,dz_{n}$.  For such $\omega$, we define $R(\omega):=\varphi_{\vert Y}$.  We can easily seen that this map is well-deined and the sequences (c) and (d) are exact.  Thus we are done.
\end{proof}

\smallskip

 \begin{notation}  We denote by $\Omega^\cdot_X((k_0+\cdot)Y)$ ($k_0$: a non-negative integer), $\Omega^{\cdot}_X(\log Y)$ and $L^\cdot(Y)$ the complexes of sheaves of $\mathbb{C}$-modules described as follows:
 \begin{equation*}
 \begin{split}
 \Omega^\cdot_X((k_0+\cdot)Y): \Omega^{0}_X(k_0Y)\to \Omega^1_X((k_0+1)Y)\to \cdots \to &\Omega^p_X((k_0+p)Y)\to \\
 &\cdots\to \Omega^n_X((k_0+n)Y),
 \end{split}
 \end{equation*}
 \begin{equation*}
 \Omega^{\cdot}_X(\log Y): \mathcal{O}_X\to \Omega^1_X(\log Y)\to \cdots\to \Omega^p_X(\log Y)\to \cdots\to 
 \Omega^n_X(\log Y),
\end{equation*}
\vskip 5pt
\hskip 15pt
 $L^\cdot(Y): \Omega^0_X\to \Phi^1_X(Y)$.
\end{notation}
\smallskip

\begin{prop}\label{prop2.2}
The natural homomorphisms of the complexes of sheaves of $\mathbb{C}$-vector spaces 
\begin{equation*}
L^\cdot(Y)\to \Omega^\cdot_X(\log Y)\to \Omega^\cdot_X((k_0+\cdot)Y)\to \Omega^\cdot_X(\ast Y)
\end{equation*}

\noindent
give rise to quasi-isomorphisms among them, and so all of the hypercohomology of these are isomorphic to $H^p(X-Y,\mathbb{C})$.
\end{prop}

\begin{proof} The former part of the proposition follows directly from Lemma \ref{lem2.1}.  The latter part is proved as follows: What we shall prove is that $\mathbb{H}^p(X,\Omega_X(\log Y))\simeq H^p(X-Y,\mathbb{C}) \quad$ ($p\geq 0$).  To do this we form a fine resolution of $\Omega^\cdot_X(\log Y)$, using 
{\it semi-meromorphic forms\/} which have poles only on $Y$.  Here, after J. Leray (\cite{Ler}), we call a $C^\infty$-differential form $\varphi$ on $X-Y$ {\it semi-meromorphic form\/} on $X$, having poles of order $k$ (at most) along $Y$ if $w^k\varphi$ is locally a $C^\infty$ regular differential form at every point of $Y$, where $w=0$ is a local defining equation of $Y$.  Similarly, as in the case of meromorphic forms, semi-meromorphic forms having {\it logarithmic poles\/} on $Y$ is defined.  We denote by $\mathfrak{A}_X^{p,q}(\log Y)$ the sheaf of germs of semi-meromorphic forms of type $(p,q)$, having {\it logarithmic poles\/} on $Y$.  Using these sheves, we obtain a fine resolution of $\Omega_X(\log Y)$ as follows:
\bigskip
\begin{equation}\label{e2.3}
\begin{CD}
\vdots & & \vdots & & \vdots & & & & \vdots \\
@AAA   @AAA  @AAA  & &  @AAA \\
\mathfrak{A}_X^{0,1}@>\text{$\partial^{0,1}$}>>\mathfrak{A}_X^{1,1}(\log Y)@>\text{$\partial^{1,1}$}>>\mathfrak{A}_X^{2,1}(\log Y)@>\text{$\partial^{2,1}$}>>\cdots @>\text{$\partial^{n,1}$}>>\mathfrak{A}_X^{n+1,1}(\log Y) \\
@AA\text{$\overline{\partial}^{0,0}$}A  @AA\text{$\overline{\partial}^{1,0}$}A 
@AA\text{$\overline{\partial}^{2,0}$}A & & 
@AA\text{$\overline{\partial}^{n+1,0}$}A \\
\mathfrak{A}_X^{0,0}@>\text{$\partial^{0,0}$}>>\mathfrak{A}_X^{1,0}(\log Y)@>\text{$\partial^{1,0}$}>>\mathfrak{A}_X^{2,0}(\log Y)@>\text{$\partial^{2,0}$}>>\cdots @>\text{$\partial^{n,0}$}>>\mathfrak{A}_X^{n+1,0}(\log Y)\\
@AAA   @AAA  @AAA  & &  @AAA \\
\mathcal{O}_X @>\text{$d$}>>\Omega_X^1(\log Y)@>\text{$d$}>>\Omega_X^2(\log Y)@>\text{$d$}>>\cdots @>\text{$d$}>> \Omega_X^{n+1}(\log Y) \\
@AAA   @AAA  @AAA  & &  @AAA \\
0 & & 0 & & 0 & & & & 0
\end{CD}
\end{equation}
\bigskip
\bigskip

\noindent
where $\mathfrak{A}_X^{p,q}$ denotes the sheave of germs of $C^\infty$ differential forms of type $(p,q)$ on $X$.  We put
\begin{eqnarray*}
A_X^{p,q}(\log Y)&:=& \Gamma(X,\mathfrak{A}_X^{p,q}(\log Y))\qquad (p\geq 0,\,\,q\geq 0),\\
A_X^{k}(\log Y)&:=& \oplus_{p+q=k}A_X^{p,q}(\log Y), \quad d^{p,q}:=\partial^{p,q}+(-1)^p\overline{\partial}^{p,q} \quad \text{and} \\
A_X^{\cdot}(\log Y)&:=& \oplus_{k}\oplus_{p+q=k}A_X^{p,q}(\log Y), \qquad d^{k}:=\oplus_{p+q=k}d^{p,q}.
\end{eqnarray*}

\noindent
Then $(A^\cdot_X(\log Y),d)$ forms a complex of $\mathbb{C}$-vector spaces and
\begin{equation*}
\mathbb{H}^p(X,\Omega_X(\log Y))\simeq H^p(A^\cdot_X(\log Y))\qquad (p\geq 0).
\end{equation*}

\noindent
By Lemma \ref{lem2.1},(d), we have the exact sequence of complexes of sheaves 
of $\mathbb{C} $-vector spaces:
\begin{equation}\label{e2.4}
0\to \Omega_X^\cdot \to \Omega_X^\cdot(\log Y)\xrightarrow{R} \Omega_Y^\cdot[-1]\to 0.
\end{equation}

\noindent
From this the following long exact sequence of hypercohomology is derived:
\begin{equation}\label{e2.5}
\to \mathbb{H}^p(\Omega_X^\cdot)\to \mathbb{H}^p(\Omega_X^\cdot(\log Y))\to 
\mathbb{H}^{p-1}(\Omega_Y^\cdot)\to \mathbb{H}^{p+1}(\Omega_X^\cdot)\to \cdots
\end{equation}

\noindent
Letting $A_X^\cdot$ and $A_Y^\cdot$ be the complexes of $\mathbb{C}$-vector spaces of global $C^\infty$ differential forms on $X$ and $Y$, respectively, we have 
$\mathbb{H}^p(\Omega_X^\cdot)\simeq H^p(A_X^\cdot)$ and $\mathbb{H}^p(\Omega_Y^\cdot)\simeq H^p(A_Y^\cdot)$.  Hence the sequence (\ref{e2.5}) is rewritten as:
\begin{equation}\label{e2.6}
\to H^p(A_X^\cdot)\xrightarrow{r^p} H^p(A_X^\cdot(\log Y))\xrightarrow{R^p} 
H^{p-1}(A_Y^\cdot)\xrightarrow{G^{p-1}} H^{p+1}(A_X^\cdot)\to \cdots
\end{equation}

\noindent
We claim that this is the dual of the homology sequence 
\begin{equation}\label{e2.7}
\leftarrow H_p(X,\mathbb{C})\xleftarrow{r_p} H_p^c(X-Y,\mathbb{C})\xleftarrow{R_{p-1}} 
H_{p-1}(Y,\mathbb{C})\xleftarrow{G_{p+1}} H_{p+1}(X,\mathbb{C})\leftarrow \cdot
\end{equation}

\noindent
(cf. (\ref{e1.3}).  In fact, since $A_X^\cdot(\log Y)$ is a subcomplex of 
$A_{X-Y}^\cdot$ which is the complex of $\mathbb{C}$-vector spaces of global 
$C^\infty$ differential forms on $X-Y$, we can define parings by integrations between the terms corresponding to each other in (\ref{e2.6}) and (\ref{e2.7}).  Furthermore, these pairings commute with the homomorphisms in 
(\ref{e2.6}) and (\ref{e2.7}), since we can easily see $A_X^\cdot(\log Y)$ is the same one as defined in Definition \ref{d1.2} and the map $R^p:H^p(A_X^\cdot(\log Y))\to H^{p-1}(A_Y^\cdot)$ is the {\it R\'esidue map\/} defined just after Definition \ref{d1.2}, and since $G^{p-1}:H^{p-1}(A_Y\cdot)\to H^{p+1}(A_X^\cdot)$ is the {\it Gysin map\/} whose description by use of differential forms has been given in Proposition \ref{prop1.7}.  Therefore, by {\it Five Lemma\/}, we conclude that the paring between $H^p(A_X^\cdot(\log Y))$ and $H^p(X-Y,\mathbb{C})$ is non-degenerated. Hence 
$H^p(A_X^\cdot(\log Y))\simeq H^p(X-Y,\mathbb{C})$.
\end{proof}
\smallskip

\begin{defn}\label{d2.1}
We define 
\begin{equation*}
I^p(X,{\ast}Y):=\Gamma(X,\Phi^p_X(\ast Y))/d\Gamma(X,\Omega_X^{q-1}(\ast Y))
\end{equation*}

\noindent
and
\begin{equation*}
I^p(X,kY):=\Gamma(X,\Phi^p_X(kY))/d\Gamma(X,\Omega_X^{q-1}((k-1)Y)).
\end{equation*}

\noindent
We call them the {\it $p$-th $\ast Y$-rational De Rham group of $X$\/} and {\it $p$-th $\ast Y$-rational De Rham group of $X$ with pole order $k$ \/}, respectively.  
\end{defn}

 Then, by Proposition \ref{prop2.2}, we have the following:
 
 \begin{prop}\label{prop2.3}
 Let $k_0$ be a positive integer such that 
 \begin{equation*}
 H^p(X,\Omega_X^q((k_0+q)Y))=0 \qquad \text{for}\quad p\geq 1,\,\,q\geq 0,
 \end{equation*}
 \end{prop}
\noindent
then,
\begin{equation*}
I^p(X,(k_0+p)Y)\simeq I^p(X,{\ast}Y)\simeq H^p(X-Y,\mathbb{C}) \qquad \text{for}\quad p\geq 0.
\end{equation*}
\smallskip

\begin{rem}\label{rem2.1}
The result in the propoition above is a special case of the theorem of 
Grothendieck (cf. \cite{Gro-1}).
\end{rem}

 Now we are going to expalin the notion of closed meromorphic forms of the {\it second kind}, having poles only along $Y$.  There are the following three different 
 definitions for this:
 
\begin{defn}\label{def2.2}
A cosed meromorphic $q$-form $\varphi$ is of the {\it second kind\/} if 
\begin{description}
\item[(A)](Picard-Lefshetz definition) at any point $x$ of $X$, there exists a 
meromorphic $q-1$ form on $X$ such that $\varphi-d\omega$ is holomorphic in a neighborhood of $x$, 
\item[(B)](Geometric R{\'e}sidue definition) it has no periods on {\it r{\'e}sidue cycles\/} (cf. Definition \ref{def1.4}) of $X-Y$, if $Y$ is sufficiently large subvariety (depending on $\varphi$), 
\item[(C)] Hodge and Atiyah's {\it algebaric definition\/}, using spectral sequences associated to the complex of sheaves of $\mathbb{C}$-vector spaces $\Omega_X^\cdot({\ast}Y)$ (or $\Omega_X^\cdot((k_0+\cdot)Y)$.
\end{description}
\end{defn}
\smallskip

We shall explain the last Hodge and Atiyah's definition (\cite{H-A}) more precisely by use of the fine resolution $\mathfrak{A}_X^{\cdot\cdot}({\ast}Y)$ of 
$\Omega_X^\cdot({\ast}Y)$, where $\mathfrak{A}_X^{\cdot\cdot}({\ast}Y)$ denotes the double complex of $\mathbb{C}$-vector spaces comprising $\mathfrak{A}_X^{p,q}({\ast}Y)$, the sheaf of germs of semi-meromorphic forms of type $(p,q)$ on $X$, having poles only along $Y$.  In the same manner as for $\mathfrak{A}_X^{p,q}(\log Y)$, we define $\mathfrak{A}_X^{p,q}(\ast Y)$ and $\mathfrak{A}_X^{k}(\ast Y)$. We form the complex of $\mathbb{C}$-vector spaces $(A_X^{\cdot}(\ast Y),d)$ for $\mathfrak{A}_X^{k}(\ast Y)$.  Then we have
\begin{equation*}
I^p(X,\ast Y)\simeq \mathbb{H}^p(X,\Omega_X^\cdot({\ast}Y))\simeq H^p(A_X^\cdot({\ast}Y)) \qquad (p\geq 0).
\end{equation*}

\noindent
under these isomorphisms, we identify $I^p(X,{\ast}Y)$ with $H^p(A_X^\cdot({\ast}Y))$ in the following.  We set
\begin{equation*}
^{\prime\prime}F^kA_X^\cdot({\ast}Y):=\oplus_{q\geq k}A_X^{\cdot q}(\ast Y),
\end{equation*}

\noindent
then $\{^{\prime\prime}F^k\}_{k\geq 0}$ give a finite decreasing filtration to $A_X^\cdot(\ast Y)$ and $A_X^\cdot(\ast Y)$ becomes a filtered complex of $\mathbb{C}$-vector spaces.  We define 

\begin{equation*}
I^p_k(X,\ast Y):=Im~\{H^p(^{\prime\prime}F^k(A_X^\cdot(\ast Y)))\to H^p(A_X^\cdot(\ast Y))\simeq I^p(X,\ast Y)~\}
\end{equation*}

\noindent
then we have a filtration on $I^p(X,\ast Y)$:
\begin{equation*}
I^p(X,\ast Y):=I^p_0(X,\ast Y)\supset I^p_1(X,\ast Y)\supset \cdots I^p_p(X,\ast Y)\supset I^p_{p+1}(X,\ast Y)=\{0\}.
\end{equation*}

\noindent
Hodge and Atiyah have defined that a closed meromorphic $p$-form $\varphi$, having poles only along $Y$, is of the {\it second kind\/} if its cohomology class 
$[\varphi]\in I^p(X,\ast Y)$ belongs to the subspace $I^p_p(X,\ast Y)$, i.e., it has the maximum filtration, and they have proved that the definitions (B) and (C) are equivalent in general. They have also proved that the definition (A) is 
equivalent to other definitions if $Y$ is a prime section of $X$.\par
\smallskip
\begin{notation} We put 
\begin{equation*}
I^p(X,\ast Y)_0=I^p_p(X,\ast Y)
\end{equation*}
\end{notation}

\noindent
Then we have:

\begin{thm}\label{thm2.4}
\hskip 5pt \par
\begin{description}
\item[(i)] $I^p(X,\ast Y)_0\simeq r^p H^p(X,\mathbb{C})\simeq H^p(X,\mathbb{C})_0 \quad (1\leseq p\leseq n+1)$,\par
\vskip 3pt
\noindent
where $r^p:H^p(X,\mathbb{C})\to H^p(X-Y,\mathbb{C})$ is the map induced by restricting closed forms on $X$ to $X-Y$, 
\item[(ii)] $I^p(X,\ast Y)_0=I^p(X,\ast Y) \quad 1\leseq p\leseq n$,
\item[(iii)] $I^n(X,\ast Y)/I^n(X,\ast Y)_0\simeq Ker~\{~H^{n-1}(Y,\mathbb{C})_0\xrightarrow{G^{n-1}} H^{n+2}(X,\mathbb{C})~\}$,\par
\vskip 3pt
\noindent
\noindent
where $G^{n-1}$ denotes the Gysin map.
\end{description}
\end{thm}

\begin{proof} Replacing $H^\ast(A_X^\cdot(\log Y))$ by $I^\ast(X,\ast Y)$ in the exact sequence (\ref{e2.6}), we obtain the exact sequence 
\begin{equation}\label{e2.8}
\to H^p(A_X^\cdot)\xrightarrow{r^p}I^p(X,\ast Y)\xrightarrow{R^p} H^{p-1}(A_Y^\cdot)\xrightarrow{G^{p-1}} H^{p+1}(A_X^\cdot)\to \cdots, 
\end{equation}

\noindent
which is dual to the homology sequence in (\ref{e2.7}).  By the {\it R{\'e}sidue definition\/} of the second kind, we have 
\begin{equation*}
I^p(X,\ast Y)_0\simeq Ann(R_{p-1}(H_{p-1}(Y,\mathbb{C}))),
\end{equation*}

\noindent
where the right hand side above denotes the annihilator subspace of $I^p(X,\ast Y)$ by $R_{p-1}(H_{p-1}(Y,\mathbb{C}))$ through the paring defined by integration between $I^p(X,\ast Y)$ and $H_p^c(X-Y,\mathbb{C})$. By the duality between 
(\ref{e2.8}) and (\ref{e2.7}), 
\begin{equation*}
 Ann(R_{p-1}(H_{p-1}(Y,\mathbb{C})))=r^{p}H^p(A_X^\cdot)\simeq r^{p}H^p(X,\mathbb{C}).
\end{equation*}

\noindent
By Proposition \ref{prop1.11}, $r^{p}H^p(X,\mathbb{C})\simeq H^p(X,\mathbb{C})_0$.  Thus we have proved (i).  By (i), (ii) follows from that $r^q:H^q(X,\mathbb{C})\to H^q(X-Y,\mathbb{C})$ is surjective for\,\, $0\leseq q\leseq n$ (cf. (\ref{e1.10}) and (\ref{e1.11})).  By the duality between (\ref{e2.8}) and (\ref{e2.7}), (iii) is trivial if we note that $R^p(I^p(X,\ast Y))\subset H^{p-1}(X-Y,\mathbb{C})_0$ (Theorem \ref{thm1.15}).
\end{proof}

\smallskip
\begin{rem}\label{rem2.2}
As in the case of $A_X^\cdot(\ast Y)$, we define a finite decreasing filtration $\{^{\prime\prime}F^k\}_{k\geq 0}$ on the complex $A_X^\cdot(\log Y)$ by 
\begin{equation*}
^{\prime\prime}F^{k}A_X^\cdot(\log Y):=\oplus_{q\geq k}A_X^{\cdot q}(\log Y)
\end{equation*}
\end{rem}
\smallskip

Then, as is well known in the homological algebra, there arises a spectral 
sequence from the filtered complex $(A_X^\cdot(\log Y),F^{\prime\prime})$ as follows:
\begin{equation*}
E_2^{p,q}:=H^p(X,\mathcal{H}^q(\Omega_X^\cdot(\log Y)))
\Longrightarrow E_\infty^{p,q}=Gr_{F^{\prime\prime}}^p\mathbb{H}^{p+q}(X,\Omega_X^\cdot(\log Y))
=Gr_{F^{\prime\prime}}^pI^{p+q}(X,\ast Y),
\end{equation*}

\noindent
where $\mathcal{H}^q(\Omega_X^\cdot(\log Y)) \quad (q\geq 0)$ are the cohomology sheaves of the complex of sheaves $\Omega_X^\cdot(\log Y)$.  From Lemma \ref{lem2.1} it follows 
\begin{equation*}
E_2^{p,q}=
\left\{
 \begin{array}{cl}
  H^p(X,\mathbb{C}) &  q=0 \\
  H^p(X,\mathbb{C}) &  q=1 \\
  0 & \text{otherwise}
  \end{array}
 \right.
\end{equation*}

\noindent
Hence we have 
\begin{equation}\label{e2.9}
\begin{split}
E_r^{q,p-q}=E_{r+1}^{q,p-q}=\cdots =E_\infty^{q,p-q}& =Gr_{F^{\prime\prime}}I^p(X,\ast Y)=0\\
 & \text{for $q\not=p, p-1$, and $r\geq 2$}
 \end{split}
\end{equation}

\noindent
This amounts to 
\begin{equation*}
 I^p(X,\ast Y)=I^p_0(X,\ast Y)=I^p_1(X,\ast Y)=\cdots =I^p_{p-1}(X,\ast Y),
\end{equation*}

\noindent
namely, the filtration of $I^p(X,\ast Y)$ induced by $\{^{\prime\prime}F^k\}_{f\geq 0}$ of $A_X^\cdot(\log Y)$ is given by a single subspace $I_p^p(X,\ast Y)$.  From this we can derive the following exact sequence (cf. \cite{God} Chapitre {\rm I}, Th{\'e}or{\`e}me 4.6.2, p.85):
\begin{equation}\label{e2.10}
\begin{CD}
\cdots@>>>E_2^{p-2,1}@>\text{$d_2^{p-1}$}>>E_2^{p,0}\\
& &  @V\text{$\simeq$}VV @V\text{$\simeq$}VV  \\
\cdots@>>> H^{p-2}(Y,\mathbb{C}) @>>> H^{p}(X,\mathbb{C}) 
 \end{CD}
\end{equation}
\begin{equation*}
\hskip 50pt\begin{CD}
\,\,\,\,\,\,@>\text{$\iota^{p}$}>>E_\infty^{p}@>\text{$j^{p}$}>>E_2^{p-1,1}@>\text{$d_2^{p}$}>>\cdots \\
& & \,\,\,\,\,\,@V\text{$\simeq$}VV @V\text{$\simeq$}VV \\
@>>>I^{p}(X,\ast Y)@>>> 
 H^{p-1}(Y,\mathbb{C})@>>>\cdots,
\end{CD}
\end{equation*}
\bigskip

\noindent
where the maps appeared in this exact sequence are described as follows:
\begin{itemize}
\item[(i)] $d_2^{p-1}$ and $d_2^p\cdots$ are the differentials of the second term $\{E_2^{p,q}\}$ of the spectral sequence, 
\item[(ii)] Since 
\begin{equation*}
\begin{split}
 E_{r+1}^{q,0}&=Ker~\{E_r^{q,0}\xrightarrow{d_r}E_r^{q+r,1-r}\}/Im~\{E_r^{q-r,r-1}\xrightarrow{d_r}E_r^{q,0}\}\\
 &= \left\{ 
 \begin{array}{rl}
 & E_r^{q,0}/Im~\{E_r^{q-r,r-1}\xrightarrow{d_r}E_r^{q,0}\} \quad r=2\\
 &  \\
 & 0 \quad r\geq 3,
 \end{array}
 \right.
 \end{split}
\end{equation*}

\noindent
there is a surjection from $E_2^{q,0}$ onto $E_\infty^{q,0}=
Gr_{^{\prime\prime}F}^qI^q(X,\ast Y)
\simeq I_q^q(X,\ast Y)$. The map $\iota^p$ is the composite of this surjection and the natural injection $I_q^q(X,\ast Y)\hookrightarrow E_\infty^q=I^q(X,\ast Y)$.
\item[(iii)] Since
\begin{equation*}
\begin{split}
 E_{r+1}^{q-1,1}&=Ker~\{E_r^{q-1,1}\xrightarrow{d_r}E_r^{q+r-1,2-r}\}/Im~\{E_r^{q-r-1,r}\xrightarrow{d_r}E_r^{q-1,1}\}\\
 &= \left\{ 
 \begin{array}{ll}
  Ker~\{E_r^{q-1,1}\xrightarrow{d_r}E_r^{q+r-1,2-r}\}& \text{\hskip 5pt $r=2$}\\ & \\
   E_r^{q-1,1} & \text{\hskip 5pt $r\geq 3$},
 \end{array}
 \right.
 \end{split}
\end{equation*}

\noindent
there is an injection from $E_\infty^{q-1,1}=Gr_{^{\prime\prime}F}^{q-1}I^q(X,\ast Y)$ into $E_2^{q-1,1}$.  The map $j^p$ is the composite of the natural surjection $E_\infty^q=I^q(X,\ast Y)$ onto $E_\infty^{q-1,1}=Gr_{^{\prime\prime}F}^{q-1}I^q(X,\ast Y)$ and the injection above from  $E_\infty^{q-1,1}=Gr_{^{\prime\prime}F}^{q-1}I^q(X,\ast Y)$ into $E_2^{q-1,1}$.  
\end{itemize}

\noindent
Chasing these maps more precisely by direct calculation, using differential forms, we can conclude that the exact sequence (\ref{e2.8}) is dual to the homology sequence (\ref{e2.7}).  Thus we have proved that $I^q(X,\ast Y)\simeq H^q(X-Y,\mathbb{C})$ again.  Besides, since the image of $\iota^p$ is $I_q^q(X,\ast Y)$ as explained above, this shows that {\it R{\'e}sidue definition\/} and {\it Hodge-Atiyah's algebraic definition\/} of the closed meromorphic forms of the second kind coincide.
\bigskip

\section{Mixed Hodge structures on $\ast Y$-rational De Rham groups of $X$ }

We call the attention of the readers to that $\Omega_X^\cdot(\log Y)$ is the most simple example of a cohomological mixed Hodge complex (CMHC) in the sense of Deligne and it induces mixed Hodge structures (MHS) on 
$\mathbb{H}^\cdot(X,\Omega_X^\cdot(\log Y))\simeq H^\cdot(X-Y,\mathbb{C})\simeq I^\cdot(X,\ast Y)$.  Concerning these MHS's a non-trivial weight filtration comes out only on $I^{n+1}(X,\ast Y)$ ($n+1={\rm dim}~X$), and it is given by a single subspace.  We shall now show that this subspace is nothing but $I^{n+1}(X,\ast Y)_0$.  First, let us recall the definition of CMHC from \cite{Elz}.  A {\it CMHC\/} $K$ on a topological space $X$ is given by 
\begin{itemize}
\item[(i)] A complex $K\in {\rm Ob}~D^{+}(X,\mathbb{Z})$ such that $H^q(X,K):=H^q(\mathbb{R}\Gamma(X,K))$ (hypercohomology of $K$) is a finite $\mathbb{Z}$-module and $H^q(X,K)\otimes \mathbb{Q}\simeq H^q(X,K\otimes \mathbb{Q})$, where $D^{+}(X,\mathbb{Z})$ denotes the derived category of lower bounded complexes of sheaves of $\mathbb{Z}$-modules over $X$.
\item[(ii)] A filtered complex $(K_{\mathbb{Q}},W)\in {\rm Ob}~D^{+}F(X,\mathbb{Q})$ and an isomorphism $K_{\mathbb{Q}}\simeq K\otimes \mathbb{Q}$ in 
$D^{+}F(X,\mathbb{Q})$ ($W$ increasing). 
\item[(iii)] A bifiltered complex $(K_{\mathbb{C}},W,F)\in {\rm Ob}~D^{+}F_2(X,\mathbb{C})$ ($W$ increasing and $F$ decreasing) and $\alpha:(K_{\mathbb{C}},W)\simeq (K_{\mathbb{Q}},W)\otimes \mathbb{C}$ in $D^{+}F(X,\mathbb{C})$, i.e., 
$Gr^W(K_{\mathbb{C}})$ and $Gr^W(K_{\mathbb{Q}})$ are quasi-isomorphic as graded comlexes, satisfying the following axioms:
\begin{itemize}
\item[(A)] $\mathbb{R}\Gamma(X,Gr_k^{W}K_{\mathbb{Q}})$, $(\mathbb{R}\Gamma(X,Gr_k^{W}K_{\mathbb{C}}),F)$ and $\mathbb{R}\Gamma(X,Gr_k^{W}\alpha):$ $\mathbb{R}\Gamma(X,Gr_k^{W}K_{\mathbb{C}})\simeq 
\mathbb{R}\Gamma(X,Gr_k^{W}K_{\mathbb{Q}})\otimes \mathbb{C}$ is a {\it Hodge complex (HC) of weight $k$\/},
\end{itemize}
\end{itemize}

\noindent
where {\it HC of weight $k$ \/} is defined as follows: A Hodge complex (HC) $K$ of weight $k$ is given by 
\begin{itemize}
\item[(i)] A complex $K\in {\rm Ob}~D^{+}(X,\mathbb{Z})$ such that the cohomology $H^q(K)$ is a $\mathbb{Z}$-module of finite type for each $q$.
\item[(ii)] A filtered complex $(K_{\mathbb{C}},F)\in {\rm Ob}~D^{+}F\mathbb{C}$ and an isomorphism $\alpha:K_{\mathbb{C}}\simeq K\otimes \mathbb{C}$ in $D^{+}\mathbb{C}$, satisfying the following axioms:
\begin{itemize}
\item[(AI)] The differential $d$ of $K_{\mathbb{C}}$ is strictly compatible to the filtration $F$, i.e., $F^i\cap Im~d=Im~(d/F^i)$ or equivalently the spectral sequence defined by $(K_{\mathbb{C}},F)$ degenerates at $E_1$ ($E_1=E_\infty$).\item[(AII)] The filtration $F$ induced on $H^q(K_{\mathbb{C}})\simeq H^q(K)\otimes \mathbb{C}$ defines a {\it HS of weight $q+k$\/}. 
\end{itemize}
\end{itemize}

\noindent
In our case, we take $K\in {\rm Ob}~D^{+}(X,\mathbb{Z})$, $(K_{\mathbb{Q}},W)\in {\rm Ob}~D^{+}F(X,\mathbb{Q})$ and $(K_{\mathbb{C}},W,F)\in {\rm Ob}~D^{+}F_2(X,\mathbb{C})$ in the definition above as follows:
\begin{equation*}
K:=\mathbb{R}j_\ast\mathbb{Z},
\end{equation*}

\noindent
where $j:X-Y \hookrightarrow X$ is the open immersion, 
\begin{eqnarray*}
& & K_{\mathbb{Q}}:=\mathbb{R}j_\ast\mathbb{Q}_{X-Y}, \\
& & W_p(K_{\mathbb{Q}}):=\tau_{\leseq p}(K_{\mathbb{Q}}),
\end{eqnarray*}

\noindent
where $\tau_{\leseq p}(K_{\mathbb{Q}})$ denotes the subcomplex of 
$K_{\mathbb{Q}}$ defined by
\begin{equation*} 
\tau_{\leseq p}(K_{\mathbb{Q}})^n=
\left\{
 \begin{array}{ll}
  K^n &  q=0 \\
  & \\
  {\rm Ker}~d &  q=1 \\
  & \\
  0 & n>p
  \end{array}
 \right.
\end{equation*}

\noindent
(which we call the {\it canonical filtration\/})
\begin{eqnarray*}
& & K_{\mathbb{C}}:=\Omega_X^\cdot(\log Y),\\
& & W_0(K_{\mathbb{C}})=\Omega_X^\cdot, \\
& & W_1(K_{\mathbb{C}})=\Omega_X^\cdot(\log Y),\\
& & F^q(K_{\mathbb{C}}):=\sigma_{\geq q}(\Omega_X^\cdot(\log Y)),
\end{eqnarray*}

\noindent
where $\sigma_{\geq q}(\Omega_X^\cdot(\log Y))$ denotes the subcomplex of 
$\Omega_X^\cdot(\log Y)$ defined by 
\begin{equation*} 
(\sigma_{\geq q}(\Omega_X^\cdot(\log Y)))^\ell=
\left\{
 \begin{array}{ll}
  0 &  \ell<q \\
  \Omega_X^{\ell}(\log Y) &  q\leseq \ell, 
  \end{array}
 \right.
\end{equation*}

\noindent
which we call the {\it stupid filtration\/}.  Instead of the filtartion $W$, we shall use the filtraion $W[q]$ defined by
\begin{equation*}
W[q]_p:=W_{p-q},
\end{equation*}

\noindent
namely, a shift by $q$ to the right on the degree of $W$.  Then $(W[q],F)$ induces a mixed Hodge structure on $H^q(\mathbb{R}\Gamma(X,\Omega_X^\cdot(\log Y)):=\mathbb{H}^q(X,\Omega_X^\cdot(\log Y))\simeq I^q(X,\ast Y)$.  We shall calculate
 $Gr_k^{W[q]}I^q(X,\ast Y)$ ($k=q, q+1$) by use of spectral sequences.  We put
\begin{eqnarray*}
& & K^\cdot:=A_X^\cdot(\log Y), \quad \text{and}\\
& & W_0(K^\cdot)=A_X^\cdot,\quad W_1(K^\cdot)=A_X^\cdot(\log Y).
\end{eqnarray*}

\noindent
$\{W_0(K^\cdot)\subset W_1(K^\cdot)=A_X^\cdot(\log Y)\}$  is the filtration induced by the filtration $\{W_0\subset W_1=\Omega_X^\cdot(\log Y)\}$ on 
$\Omega_X^\cdot(\log Y)$.  We define 
\begin{equation*}
W_p^\prime(K^\cdot):=W[q]_{-p}(K^\cdot)=W_{-p-q}(K^\cdot)\quad (p\leseq -q).
\end{equation*}

\noindent
Then $\{W_p^\prime(K^\cdot)\}$ is a decreasing filtration of $K^\cdot$.  Hence we can consider the spectral sequence concerning the filtration complex $(K^\cdot,W^\prime(K^\cdot))$, whose 
$0$-th term and $1$-st one are computed as follows:
\begin{eqnarray*}
_{W^\prime} E_0^{r,s}&=&Gr^r_{W^\prime}(K^{r+s})\\
     &=& \left\{
        \begin{array}{ll}
        W_0(K^{s-q}) &  r=-q \\
        & \\
        W_1(K^{s-q-1})/W_0(K^{s-q-1}) &  r=-q-1 \\
        & \\
        0 & \text{otherwise}
        \end{array}
        \right.\\
     &=&  \left\{
         \begin{array}{ll}
         A_X^{s-q} & r=-q \\
         & \\
         A_X^{s-q-1}(\log Y)/A_X^{s-q-1}\simeq A_Y^{s-q-2} &  r=-q-1\\
         & \\
         0 & \text{otherwise},
         \end{array}
         \right.
\end{eqnarray*}

\noindent
where the isomorphism $A_X^{s-q-1}(\log Y)/A_X^{s-q-1}\simeq A_Y^{s-q-2}$ comes from the exact sequence of sheaves 
\begin{equation*}
0\to \mathcal{W}_X^\cdot \to \mathcal{A}_X^\cdot(\log Y)\xrightarrow{R}\mathcal{A}_Y[-1]^\cdot\to 0,
\end{equation*}

\noindent
(cf. Proposition \ref{prop1.7}) which is the $C^\infty$ version of the exact sequence (\ref{e2.4});
\begin{eqnarray*}
_{W^{\prime}}E_1^{r,s}&=&\frac{{\rm Ker}\{_{W^{\prime}}E_1^{r,s}\xrightarrow{d_1}_{W^{\prime}}E_1^{r+1,s}\}}{{\rm Im}\{_{W^{\prime}}E_1^{r-1,s}\xrightarrow{d_1}_{W^{\prime}}E_1^{r,s}\}}\\
 & & \\
 &=& \left\{
 \begin{array}{ll}
 &  H^{s-q}(A_X^\cdot) \simeq H^{s-q}(X,\mathbb{C}_X) \quad r=-q,\,\,s\geq q\\
 & \\
 & H^{s-q-1}(A_X^\cdot(\log Y)/A_X^\cdot)\simeq 
    H^{s-q-1}(A_X^\cdot(\log Y)/W_X^\cdot)\\ 
  & \simeq H^{s-q-2}(A_Y^\cdot)\simeq H^{s-q-2}(Y,\mathbb{C}_Y) 
  \quad r=-q-1,\,\, s\geq q+1, \\
& \\
 & 0 \qquad \text{otherwise}.
 \end{array}
 \right.
\end{eqnarray*}

\noindent
Hence we have 
\begin{equation}\label{e3.1}
\begin{split}
 _{W^{\prime}}E_t^{r,p-r}=_{W^{\prime}}E_{t+1}^{r,p-r}=\cdots & =_{W^{\prime}}E_\infty^{r,p-r}=Gr_r^{W^\prime}I^p(X,\ast Y)=0\\
 & \text{for}\quad r\not=-q, -q-1\quad \text{and}\quad t\geq 2
 \end{split}
\end{equation}

\noindent
This is equivalent to $I^p(X,\ast Y)=W^\prime_{-q-1}(I^p(X,\ast Y))\supset 
W_{-q}^\prime(I^p(X,\ast Y))=E_\infty^{-q,p+q}$. From these we obtain the following exact sequnece:

\begin{equation}\label{e3.2}
\begin{CD}
\cdots @>>> _{W^\prime}E_1^{-q,p}@>\text{$\iota^p$}>>_{W^\prime}E_\infty^{p,-q}
@>\text{$\j^p$}>>_{W^\prime}E_1^{-q-1,p+1}\\
& &  @V\text{$\simeq$}VV  @V\text{$\simeq$}VV  @V\text{$\simeq$}VV \\
\cdots@>>>H^{p-q}(X,\mathbb{C})@>\text{$r^{p-q}$}>> I^{p-q}(X,\ast Y)@>\text{$R^{p-q}$}>>H^{p-q-1}(Y,\mathbb{C})
\end{CD}
\end{equation}
\begin{equation*}
\text{\hskip 50pt} \begin{CD}
\,\,\,\,\,\,\,\,\,@>\text{$d_1^{p+1}$}>>_{W^\prime}E_1^{-q,p+1}@>\text{$\iota^{p+1}$}>>_{W^\prime}E_\infty^{p,-q+1}@>>>\cdots \\
\,\,\,\,\,\,\,\,\, & &     @V\text{$\simeq$}VV  @V\text{$\simeq$}VV & & \\
\,\,\,\,\,\,\,\,\, @>\text{$G^{p-q-1}$}>>H^{p-q+1}(X,\mathbb{C})@>\text{$r^{p-q+1}$}>>I^{p-q+1}(X,\ast Y)@>>>\cdots,
\end{CD}
\end{equation*}  

\noindent
where the maps in this diagram are described as follows:
\begin{itemize}
\item[(i)] $d_1^{p+1}$ is the differential at the first term of $\{E_1^{p-q+1,p+1}\}$ of the spectral sequnece,
\item[(ii)] Since 
\begin{eqnarray*}
_{W^{\prime}}E_{r+1}^{-q,p}&=&\frac{{\rm Ker}\{_{W^{\prime}}E_r^{-q,p}\xrightarrow{d_r}_{W^{\prime}}E_r^{-q+r,p-r+1}\}}{{\rm Im}\{_{W^{\prime}}E_r^{-q-r,p+r-1}\xrightarrow{d_r}_{W^{\prime}}E_r^{-q,p}\}}\\
 & & \\
 &=& \left\{
 \begin{array}{ll}
 & _{W^{\prime}}E_r^{-q,p}/{\rm Im}\{_{W^{\prime}}E_r^{-q-1,p}\xrightarrow{d_r}_{W^{\prime}}E_r^{-q,p}\}, \quad r=1\\
 & \\
 & _{W^{\prime}}E_r^{-q,p}, \quad r\geq 2,
 \end{array}
 \right.
\end{eqnarray*}

\noindent
there is a surjection from $_{W^{\prime}}E_1^{-q,p}$ onto $_{W^{\prime}}E_\infty^{-q,p}=Gr_{W^{\prime}}^{-q}I^{p-q}(X,\ast Y)=W_{-q}^{\prime}I^{p-q}(X,\ast Y)$.  The map $\iota^p$ is the composite of this surjection and the natural injection $W_{-q}^{\prime}I^{p-q}(X,\ast Y)\hookrightarrow I^{p-q}(X,\ast Y)=_{W^{\prime}}E_\infty^{-q,p}$.
\item[(iii)] Since 
\begin{eqnarray*}
_{W^{\prime}}E_{r+1}^{-q-1,p+1}&=&\frac{{\rm Ker}\{_{W^{\prime}}E_r^{-q-1,p+1}\xrightarrow{d_r}_{W^{\prime}}E_r^{-q-1+r,p-r+2}\}}{{\rm Im}\{_{W^{\prime}}E_r^{-q-1-r,p+r}\xrightarrow{d_r}_{W^{\prime}}E_r^{-q-1,p+1}\}}\\
 & & \\
 &=& \left\{
 \begin{array}{ll}
 & {\rm Ker}\{_{W^{\prime}}E_1^{-q-1,p+1}\xrightarrow{d_1}_{W^{\prime}}E_1^{-q,p+1}\}, \quad r=1\\
 & \\
 & _{W^{\prime}}E_r^{-q-1,p+1}, \quad r\geq 2,
 \end{array}
 \right.
\end{eqnarray*}

\noindent
there is an injection from $_{W^{\prime}}E_\infty^{-q-1,p+1}=Gr_{W^\prime}^{-q-1}I^{p-q}(X,\ast Y)$ into $_{W^{\prime}}E_1^{-q-1,p+1}$.  The map $j^p$ is the composite of the natural surjection $_{W^{\prime}}E_\infty^{p-q}=I^{p-q}(X,\ast Y)$ onto $_{W^{\prime}}E_\infty^{-q-1,p+1}Gr_{W^\prime}I^{p-q}(X,\ast Y)$ and the injection above from $_{W^{\prime}}E_\infty^{-q-1,p+1}$ into $_{W^{\prime}}E_1^{-q-1,p+1}$. 
\end{itemize}

\noindent
Chasing these maps more precisely by direct calculation, using differential forms, we can conclude that the exact sequence (\ref{e3.2}) is dual to the homology sequence (\ref{e2.7}).  By the definition of the map $\iota^p$ and $j^p$, 
we have 
\begin{eqnarray*}
\iota^p(_{W^{\prime}} E_1^{-q,p})&=& W^\prime_{-q}I^{p-q}(X,\ast Y) \quad \text{and}\\
j^p(_{W^{\prime}} E_\infty^{p+q})&=& Gr_{W^\prime}^{-q-1}I^{p-q}(X,\ast Y),
\end{eqnarray*}

\noindent
which are rewritten as 
\begin{eqnarray*}
r^{p-q}(H^{p-q}(X,\mathbb{C})&=& W^\prime_{-q}I^{p-q}(X,\ast Y) \quad \text{and}\\
R^{p-q}(I^{p-q}(X,\ast Y))&=& Gr_{W^\prime}^{-q-1}I^{p-q}(X,\ast Y),
\end{eqnarray*}

\noindent
If we put $p=2q$, then we have 
\begin{equation*}
H^q(X,\mathbb{C})_0=r^q(H^q(X,\mathbb{C}))=W_{-q}^{\prime}I^q(X,\ast Y)=
W[q]_qI^q(X,\ast Y)
\end{equation*}

\noindent
and
\begin{eqnarray*}
{\rm Ker}\{G^{q-1}:H^{q-1}(Y,\mathbb{C})_0\to F^kH^{q+2}(X,\mathbb{C})\}&=&
R^q(I^q(X,\ast Y)) \\
&\simeq& Gr_{W^\prime}^{-q-1}I^{q}(X,\ast Y) \\
&\simeq& =Gr^{W[q]}_{q+1}I^{q}(X,\ast Y).
\end{eqnarray*}

\noindent
Therefore, combining these results with those of Theorem \ref{thm2.4}, we have
\begin{thm}\label{thm3.1}
\hskip 50pt
\begin{itemize}
\item[{\rm (i)}] $Gr_q^{W[q]}H^q(X-Y,\mathbb{C})=W[q]_qH^q(X-Y,\mathbb{C})=I^q(X,\ast Y)_0$,
\item[{\rm (ii)}] $Gr_{q^+1}^{W[q]}H^q(X-Y,\mathbb{C})=I^q(X,\ast Y)/I^q(X,\ast Y)_0$,
\item[{\rm (iii)}] $F^kGr_{q}^{W[q]}H^q(X-Y,\mathbb{C})\simeq F^kH^q(X,\mathbb{C})_0$,
\item[{\rm (iv)}] $F^kGr_{q^+1}^{W[q]}H^q(X-Y,\mathbb{C})
\simeq{\rm Ker}\{F[-1]^kH^{q-1}(Y,\mathbb{C})_0\xrightarrow{G^{q-1}}F^kH^{q+2}(Y,\mathbb{C})\}$
\end{itemize}
\end{thm}

 From now on, we consider the following complex of sheaves of $\mathbb{C}$-vector spaces:
\begin{equation*}
\Omega_X^\cdot((1+\cdot)Y): \mathcal{O}_X(Y)\to \Omega_X^1(2Y)\to \cdots \to 
\Omega_X^p((p+1)Y)\to \cdots \to \Omega_X^{n+1}((n+2)Y)
\end{equation*}

\noindent
We define a decreasing filtration $\{F^{\prime k}\}_{0\leseq k\leseq n}$ by 
\begin{equation}\label{e3.3}
\begin{split}
F^{\prime k}(\Omega_X^\cdot((1+\cdot)Y)):=\{\cdots & \to 0\to \Omega_X^k(Y)\to 
\Omega_X^{k+1}(2Y)\to \\
 & \cdots  \to \Omega_X^p((p-k+1)Y)\to\cdots \to  \Omega_X^{n+1}((n-k+2)Y)\},
\end{split}
\end{equation}

\noindent
and an increasing filtartion $\{W_0^\prime\subset W_1^\prime\}$ by 
\begin{eqnarray*}
W_0^\prime(\Omega_X^\cdot((1+\cdot)Y))&:& \mathcal{O}_X\to 
\Omega_X^{1}\to\cdots \to \Omega_X^p\to\cdots \to \Omega_X^{n+1},\\
& & \\
W_1^\prime(\Omega_X^\cdot(1+\cdot))&:& \Omega_X^\cdot(1+\cdot).
\end{eqnarray*}

\noindent
Then we have

\begin{prop}\label{prop3.2}
The bi-filtered complexes of sheaves of $\mathbb{C}$-vector spaces 
\begin{equation*}
(\Omega_X^\cdot(\log Y),W,F)\quad \text{and}\quad (\Omega_X^\cdot((1+\cdot)Y),W^\prime,F^\prime)
\end{equation*}
\noindent
are quasi-isomorphic, i.e., bi-graded complex of sheaves of $\mathbb{C}$-vector spaces 
\begin{equation*}
Gr_{F}Gr^{W}(\Omega_X^\cdot(\log Y))\quad \text{and}\quad 
Gr_{F^\prime}Gr^{W^\prime}(\Omega_X^\cdot((1+\cdot)Y))
\end{equation*}
\noindent
are quasi-isomorphic where the filtration $F$ of 
$\Omega_X^\cdot(\log Y)$ is defined by 
\begin{equation*}
\begin{split}
F^{k}(\Omega_X^\cdot(\log Y)):=\{\cdots\to 0\to \Omega_X^k(\log Y)\to 
\Omega_X^{k+1}(\log Y)\to \cdots \to & \Omega_X^{n+1}(\log Y)\}\\
& (0\leseq k\leseq n+1),
\end{split}
\end{equation*}
\end{prop}

\begin{proof} First, we have 
\begin{eqnarray*}
Gr_{F^\prime}^kGr_0^{W^\prime}(\Omega_X^\cdot((1+\cdot)Y))&=& \Omega_X^k[-k], \\ & & \\
Gr_{F^\prime}^kGr_1^{W^\prime}(\Omega_X^\cdot((1+\cdot)Y))&=& (\Omega_X^{k+\cdot}((1+\cdot)Y)/\Omega_X^{k+\cdot}(\cdot Y))[-k], \\
& & \\
Gr_F^kGr_0^{W}(\Omega_X^\cdot(\log Y))&=& \Omega_X^k[-k], \quad \text{and}\\
& & \\
Gr_F^kGr_1^{W}(\Omega_X^\cdot(\log Y))&=&
(\Omega_X^{k+\cdot}(\log Y)/\Omega_X^{k+\cdot})[-k]  \\
\end{eqnarray*}

\noindent
Thus $Gr_{F^\prime}^kGr_0^{W^\prime}(\Omega_X^\cdot((1+\cdot)Y))$ and $Gr_F^kGr_0^{W}(\Omega_X^\cdot(\log Y))$ are quasi-isomorphic and
\begin{equation*}
H^p(Gr_F^kGr_1^{W}(\Omega_X^\cdot(\log Y))=
\left\{
\begin{array}{cl}
\Omega_X^p(\log Y)/\Omega_X^p & p=k\geq 1\\
0 & \text{otherwise}.
\end{array}
\right.
\end{equation*}

\noindent
We shall calculate $H^p(Gr_{F^\prime}^kGr_1^{W^\prime}
(\Omega_X^\cdot((1+\cdot)Y))$.  Obviously, 
\begin{equation*}
H^p(Gr_F^kGr_1^{W^\prime}(\Omega_X^\cdot((1+\cdot)Y))=0\quad \text{for}\quad 0\leseq p\leseq k-1, \quad 1\leseq k. 
\end{equation*}

\noindent
Assume $p\geq k+1$.  Let 
$[\omega]\in \Omega_X^p((p-k+1)Y)/\Omega_X^p((p-k)Y)$ 
be an element with $d[\omega]=0$ in $\Omega_X^p((p-k+2)Y)/\Omega_X^p((p-k+1)Y)$ where $\omega$ is an element of $\Omega_X^{p+1}((p-k+1)Y)$.  Since $d\omega$ is a closed form, by Lemma \ref{lem2.1}, (i)-(a), there exists $\varphi\in \Omega_X^p((p-k)Y)$ such that $d\varphi=d\omega$.  Since $\omega-\varphi\in \Phi_X^p((p-k+1)Y)$, by the same reason, there exists $\psi\in \Omega_X^{p-1}((p-k)Y)$ such 
that $d\psi=\omega-\varphi$.  This means $d[\psi]=[\omega]$.  Thus 
$H^p(Gr_F^kGr_1^{W^\prime}(\Omega_X^\cdot((1+\cdot)Y))=0$ for $p\geq k+1$.  Let 
$[\omega]\in \Omega_X^k(Y)/\Omega_X^{k}$ be an element with $d[\omega]=0$ in 
$\Omega_X^{k+1}(2Y)/\Omega_X^{k+1}(Y)$.  This amounts to $d\omega\in \Omega_X^{k+1}(Y)$.  If $k\geq 1$, we can easily see that this is the case if and only if $\omega\in \Omega_X^{k}(\log Y)$.  This fact tells us that
\begin{equation*}
H^k(Gr_{F^\prime}^kGr_1^{W^\prime}(\Omega_X^\cdot((1+\cdot)Y))\simeq \Omega_X^k(\log Y)/\Omega_X^k \quad \text{for}\quad k\geq 1.
\end{equation*}

\noindent
If $k=0$, we can easily see that $\omega\in \mathcal{O}_X$, since  
$\omega\in \mathcal{O}_X(Y)$, $d\omega\in \Omega_X^1(Y)$.  Hence 
$H^0(Gr_{F^\prime}^0Gr_1^{W^\prime}(\Omega_X^\cdot(1+\cdot)))=0$. 
This completes the proof.
\end{proof}
\smallskip

 We define
 \begin{equation*}
 I^p_k(X,(p+1)Y):=\frac{\Gamma(X,\Phi_X^p((p-k+1)Y))}{d\Gamma(X,\Omega_X^{p-1}((p-k)Y))} \quad (0\leseq k\leseq p)
 \end{equation*}

\noindent
and denote by $I_k^p(X,(p+1)Y)_0$ the subspace of $I_k^p(X,(p+1)Y)$ generated 
by closed moromorphic of $p$-forms of the second kind.  The CMHC 
$(\Omega_X^\cdot(\log Y),W,F)$ induces a mixed Hodge structure on $H^p(X-Y,\mathbb{C})$ ($\simeq \mathbb{H}^p(X,\Omega_X(\log Y))$).  We denote by 
$\{F^kH^p(X-Y,\mathbb{C})\}_{0\leseq k\leseq p}$ the Hodge filtration of $H^p(X-Y,\mathbb{C})$ concerning this mixed Hodge structure, and by $\{F^kH^p(X,\mathbb{C})_0\}_{0\leseq k\leseq p}$ the ordinary Hodge filtration of $H^p(X,\mathbb{C})_0$, the $p$-th primitive cohomology group of $X$.  With this notation we have

\begin{thm}\label{thm3.3}
If $Y$ is sufficiently ample so that 
\begin{equation}\label{e3.4}
H^p(X,\Omega_X^q(kY))=0 \quad \text{for} \quad p\geq 1, q\geq 0, k\geq 1,
\end{equation}

\noindent
 then we have
\begin{equation}\label{e3.5}
F^kH^p(X-Y,\mathbb{C}) \simeq I^p_k(X,(p+1)Y) \quad 0\leseq k\leseq p \quad \text{and}
\end{equation}
\begin{equation} \label{e3.6}
F^kH^p(X,\mathbb{C})_0 \simeq  I^p_k(X,(p+1)Y)_0 \quad 0\leseq k\leseq p 
\end{equation}

\noindent
under the isomorphisms $H^p(X-Y,\mathbb{C})\simeq I^p(X,(p+1)Y)$ and 
$H^p(X,\mathbb{C})_0\simeq I^p(X,(p+1)Y)_0$ in Proposition \ref{prop2.3} and 
Theorem \ref{thm2.4}, respectively. 
\end{thm}

\begin{proof} Using the sheaves $\mathfrak{A}_X^{p,q}(\ell Y)$, the sheaves of germs of semi-meromorphic forms of type $(p,q)$ on $X$, having poles of order $\ell$ (at most) alomg $Y$, we can form a fine resolution of it by use of more small 
sheaves.  Let $\mathfrak{B}_X^{p,q}(\ell Y)$ be the subsheaves of 
$\mathfrak{A}_X^{p,q}(\ell Y)$ characterized by the following prescription: Letting $\varphi$ be a local cross-section of $\mathfrak{A}_X^{p,q}(\ell Y)$ if 
and only if $f^{\ell-1}df\wedge \varphi$ is a $C^\infty$ regular differential form where $f=0$ is a local holomorphic defining equation for $Y$.  Using $\mathfrak{B}_X^{p,q}(\ell Y)$, we obtain a fine resolution of $\Omega_X^\cdot((1+\cdot)Y)$ as follows:
\bigskip
\begin{equation}\label{e3.7}
\begin{CD}
\vdots & & \vdots & & \vdots & & & & \vdots \\
@AAA   @AAA  @AAA  & &  @AAA \\
\mathfrak{B}_X^{0,1}(Y)@>\text{$\partial^{0,1}$}>>\mathfrak{B}_X^{1,1}(2Y)@>\text{$\partial^{1,1}$}>>\mathfrak{B}_X^{2,1}(3Y)@>\text{$\partial^{2,1}$}>>\cdots @>\text{$\partial^{n,1}$}>>\mathfrak{B}_X^{n+1,1}((n+2)Y) \\
@AA\text{$\overline{\partial}^{0,0}$}A  @AA\text{$\overline{\partial}^{1,0}$}A 
@AA\text{$\overline{\partial}^{2,0}$}A & & 
@AA\text{$\overline{\partial}^{n,0}$}A \\
\mathfrak{B}_X^{0,0}(Y)@>\text{$\partial^{0,0}$}>>\mathfrak{B}_X^{1,0}(2Y)@>\text{$\partial^{1,0}$}>>\mathfrak{B}_X^{2,0}(3Y)@>\text{$\partial^{2,0}$}>>\cdots @>\text{$\partial^{n,0}$}>>\mathfrak{B}_X^{n+1,0}((n+2)Y)\\
@AAA   @AAA  @AAA  & &  @AAA \\
\Omega_X^0(Y) @>\text{$d$}>>\Omega_X^1(2Y)@>\text{$d$}>>\Omega_X^2(3Y)@>\text{$d$}>>\cdots @>\text{$d$}>> \Omega_X^{n+1}((n+2)Y) \\
@AAA   @AAA  @AAA  & &  @AAA \\
0 & & 0 & & 0 & & & & 0
\end{CD}
\end{equation}
\vskip 10pt

\noindent
We put 
\begin{eqnarray*}
B_X^{p,q}((p+1)Y)&:=& \Gamma(X,\mathfrak{B}_X^{p,q}((p+1)Y)\quad (p\geq 0, q\geq 0), \\
B_X^k((k+1)Y)&:=& \oplus_{p+q=k}B_X^{p,q}((p+1)Y)\quad d^{p,q}=\partial^{p,q}+(-1)^p\overline{\partial}^{p,q}\quad \text{and}\\
B_X^\cdot((1+\cdot)Y)&:=& \oplus_k\oplus_{p+q=k}B_X^{p,q}((p+1)Y)
\end{eqnarray*}

\noindent
Then $(B_X^\cdot((1+\cdot)Y),d)$ forms a complex of $\mathbb{C}$-vector 
spaces and we have 
\begin{equation*}
\mathbb{H}^p(X,\Omega_X((1+\cdot)Y))\simeq H^p(B_X^\cdot((1+\cdot)Y)) \quad 
(p\geq 0).
\end{equation*}
 
\noindent
The filtration $\{F^{\prime k}\}$ of $\Omega_X((1+\cdot)Y)$ defined in (\ref{e3.3}) induces a filtration on $B_X^\cdot((1+\cdot)Y)$, which we denote by 
$\{F^{\prime k}B_X^\cdot((1+\cdot)Y)\}$, i.e.,
\begin{equation*}
F^{\prime k}B_X^\cdot((1+\cdot)Y):=\oplus_p\oplus_{p\geq q\geq k}B_X^{q,p-q}((q+1)Y)
\end{equation*}

\noindent
Since $(\Omega_X((1+\cdot)Y),W^\cdot,F^\prime)$ is a CMHC by Proposition \ref{prop3.2}, the spectral sequence, associated to the filtration $\{F^{\prime k}B_X^\cdot((1+\cdot)Y)\}$ and whose final terms are 
\begin{equation*}
_{F^\prime}E_\infty^{p,q}=Gr^p_{F^\prime}=Gr_{F^\prime}^pH^{p+q}(B_X^{\cdot}((1+\cdot)Y)),
\end{equation*}

\noindent
is degenerated at the $1$-st term (cf. \cite{De:II}, Th{\'e}or{\`e}me 3.2.5, \cite{Elz}, Th{\'e}or{\`e}me 3.2.1).  Therefore, we have 
\begin{equation}\label{e3.8}
\begin{split}
_{F^\prime}E_1^{k,p-k}&=H^p(F^k(B^\cdot)/F^{k+1}(B^\cdot)) \quad (B^\cdot=
B_X^\cdot((1+\cdot)Y))\\
    & \simeq _{F^\prime}E_\infty^{k,p-k}=Gr_{F^\prime}^kH^p(B^\cdot)
\end{split}
\end{equation}

\noindent
Here we should recall that the filtration on $H^p(B^\cdot)$ induced by 
$\{F^\prime\}$ on $B^\cdot$ is defined by 
\begin{eqnarray*}
F^{\prime k}H^p(B^\cdot)&:=&{\rm Im}\{H^p(F^k(B^\cdot))\to H^p(B^\cdot)\} \quad \text{and}\\
Gr_{F^{\prime}}^kH^p(B^\cdot)&=& F^{\prime k}H^p(B^\cdot)/F^{\prime k+1}H^p(B^\cdot)
\end{eqnarray*}

\noindent
From this and (\ref{e3.8}) it follows that the natural map
\begin{equation*}
H^p(F^{\prime k}(B^\cdot))\to H^p(F^{\prime k}(B^\cdot)/F^{\prime k+1}(B^\cdot))
\end{equation*}

\noindent
is surjective.  Hence the long exact sequence of cohomology associated to 
the exact sequence of complex 
\begin{equation*}
0\to F^{\prime k+1}(B^\cdot)\to F^{\prime k}(B^\cdot)\to 
F^{\prime k}(B^\cdot)/F^{\prime k+1}(B^\cdot)\to 0
\end{equation*}

\noindent
breaks up into the following short exact sequences 
\begin{equation*}
\begin{split}
0 \to H^p(F^{\prime k+1}(B^\cdot))\to H^p(F^{\prime k}(B^\cdot))\to H^p(F^{\prime k}(B^\cdot)/F^{\prime k+1}(B^\cdot))& \to 0 \\
&(0\leseq p\leseq n+1, 0\leseq k\leseq p)
\end{split}
\end{equation*}

\noindent
Here $H^p(F^{\prime k}(B^\cdot)/F^{\prime k+1}(B^\cdot))\simeq  Gr_{F^\prime}^kH^p(B^\cdot)$.  Hence 
\begin{equation}\label{e3.9}
H^p(F^{\prime k}(B^\cdot)\simeq F^{\prime k}H^p(B^\cdot)\simeq F^{\prime k}H^p(X-Y,\mathbb{C}) \quad (0\leseq k\leseq p, 0\leseq p\leseq n+1). 
\end{equation}

\noindent
On the other hand, by the assumption \ref{e3.4}, we have
\begin{equation}\label{e3.10}
\begin{split}
H^p(F^{\prime k}(B^\cdot))&\simeq \mathbb{H}^p(F^{\prime k}(\Omega_X^{\cdot}((1+\cdot)Y)))\\
 &=\mathbb{H}^p(\Omega_X^{k+\cdot}((1+\cdot)Y)[-k]) \\
 &=\mathbb{H}^{p-k}(\Omega_X^{k+\cdot}((1+\cdot)Y))\\
 &\simeq \frac{\Gamma(X,\Phi_X^p((p-k+1)Y))}{d\Gamma(X,\Omega_X^{p-1}((p-k)Y))}
\\
&=I_k^p(X,(p+1)Y).
\end{split}
\end{equation}

\noindent
By Proposition \ref{prop3.2}, the ordinary Hodge filtartion $F^kH^p(X-Y,\mathbb{C})$ of the cohomology $H^p(X-Y,\mathbb{C})$ coincides with $F^{\prime k}H^p(X-Y,\mathbb{C})$.  Therefore, by (\ref{e3.9}) and (\ref{e3.10}), we conclude that (\ref{e3.5}) certainly holds.  Noticing that $I^p(X,(p+1)Y)_0\simeq I^p(X,\ast Y)_0$, 
we obtain (\ref{e3.6}) from (\ref{e3.5}) and Theorem \ref{thm3.1}.
\end{proof}
\bigskip

\section{Generalized Poincar{\'e} r{\'e}sidue map }
 
  The setting under which we shall work in this section is as follows:
  Let $X$ be a non-singular irreducible algebraic variety of dimension 
 $n+1$ embedded in a sufficiently higher complex projective space 
 $\mathbb{P}^{N}$, $Y$ a generic hyperplane section of $X$ which satisfies 
 the condition (\ref{e3.4}) in Theorem \ref{thm3.3}, and $Y^\prime$ a non-singular, irreducible hypersurface section of sufficiently higher degree such that if we set $Z=Y\cdot Y^\prime$, then 
\begin{equation}\label{e4.1}
H^p(Y,\Omega_Y^q(kZ))=0 \quad \text{for} \quad p\geq 1, q\geq 0 \quad \text{and}\quad 
k\geq 1.
\end{equation}

\noindent
When we refer to primitive cohomology, we always means the one concerning the Hodge metric whose fundamental forms is dual to the homology class $[Y]$ (resp. $[Z]$).  Under this setting and with the same notation as in the previous sections, the purpose of this section is to define the so-called {\it generalized Poinca{\'e} residue map\/}
\begin{equation*}
R\acute{e}s: I^{n+1}(X,(n+2)Y)\to I^{n}(Y,(n+1)Z)_0
\end{equation*}

\noindent
and prove the following theorem:
\begin{thm}\label{thm4.1}
Under the setting above, we have 
\begin{eqnarray*}
F^kH^n(Y,\mathbb{C})_0&\simeq& I^n_k(Y,(n+1)Z)_0 \\
                      &\simeq& R\acute{e}s(I^{n+1}_{k+1}(X,(n+2)Y))\oplus r^n(I^n_k(X,(n+1)Y^\prime)_0)),
\end{eqnarray*}
\end{thm}

\noindent
where $r^n$ denote the map induced by the natural map 
$H^n(X,\mathbb{C})_0\to H^n(Y,\mathbb{C})_0$. \par
\bigskip

We shall prove the theorem after several lemmas and 
Propositions. We denote by $\Omega^q(kY+\ast Y^\prime)$ the sheaf of germs of 
meromorphic $q$-forms having poles of order $k$ (at most) along $Y$ and poles of arbitrary order along $Y^\prime$ as their only singularities.  We denote by $\Omega^q(\log Y+kY^\prime)$ the sheaf of germs of meromorphic $q$-forms having logarithmic poles along $Y$ and poles of order $k$ at most as their only singulatities.
We consider the following homomorphisms of comlexes of sheaves of $\mathbb{C}$-vector spaces:
\begin{equation}\label{e4.2}
\begin{CD}
 \Omega_X^\cdot((1+\cdot)Y): & & \mathcal{O}_X(Y)@>>>\Omega_X^1(2Y)@>>>\cdots \\
& & @VVV  @VVV & & \\
\Omega_X^\cdot((1+\cdot)Y+\ast Y^\prime)): & & \mathcal{O}_X(Y+\ast Y^\prime)@>>>\Omega_X^1(2Y+\ast Y^\prime)@>>>\cdots \\
& & @AAA  @AAA & & \\
\Omega_X^\cdot(\log Y+(1+\cdot)Y^\prime): & & \mathcal{O}_X(Y^\prime)@>>>\Omega_X^1(\log Y+2Y^\prime)@>>>\cdots
\end{CD}
\end{equation}
\begin{equation*}
\hskip 5pt \begin{CD}
@>>>\Omega_X^p((p+1)Y)@>>>\cdots@>>>\Omega_X^{n+1}((n+2)Y) \\
 & & @VVV & & @VVV \\
@>>>\Omega_X^p((p+1)Y+\ast Y^\prime)@>>>\cdots@>>>\Omega_X^{n+1}((n+2)Y+\ast Y^\prime) \\ & & @AAA & & @AAA \\
@>>>\Omega_X^p(\log Y+(p+1)Y^\prime)@>>>\cdots@>>>
\Omega_X^{n+1}(\log Y+(n+2)Y^\prime)
\end{CD}
\end{equation*}
\smallskip

\begin{prop}\label{prop4.2}
The homomorphism of complexes of sheaves 
\begin{equation*}
\Omega_X^\cdot(\log Y+(1+\cdot)Y^\prime)\to \Omega_X^\cdot((1+\cdot)Y+\ast Y^\prime)
\end{equation*}
 in the diagram (\ref{e4.2}) is a quasi-isomorphism.
\end{prop}
\begin{proof} By virtue of Proposition \ref{prop2.2} it suffices to show that the stalks of the cohomology sheves $\mathcal{H}^p(\Omega_X^\cdot(\log Y+(1+\cdot)Y^\prime))$ and $\mathcal{H}^p(\Omega_X^\cdot((1+\cdot)Y+\ast Y^\prime)$ are isomorphic at a point $x_0\in Y\cap Y^\prime$.  Let $(z_1,\cdots,z_{n+1})$ be a holomorphic local coordinate system at $x_0$ such that $z_1=0$ and $z_2=0$ are local defining equations $Y$ and $Y^\prime$, respectively.  We are going to show that 
\begin{equation}\label{e4.3}
\begin{split}
&\mathcal{H}^p(\Omega_X^\cdot((1+\cdot)Y+\ast Y^\prime)) \\
& \simeq \mathcal{H}^p(\Omega_X^\cdot(\log Y+(1+\cdot)Y^\prime))=
\left\{
\begin{array}{cl}
\mathbb{C}_X &  p=0 \\
 & \\
\mathbb{C}\{\displaystyle \frac{dz_1}{z_1},\frac{dz_2}{z_2}\} & p=1 \\
 & \\
\mathbb{C}\{\displaystyle \frac{dz_1dz_2}{z_1z_2}\} & p=2 \\
& \\
0 & \text{otherwize}
\end{array}
\right.
\end{split}
\end{equation}

\noindent
Now let $\varphi=dz_1\wedge \alpha+\beta$ be a local cross-section of 
$\Phi_X^p((p+1)Y+\ast Y^\prime)$ \quad ($p\geq 1$) in a neighborhood of 
$x_0$, where $\alpha,\beta$ are local meromorphic forms, having poles of order 
$p+1$ (at most) along $Y$ and poles of arbitrary order along $Y^\prime$ as their only singularities, and not involving $dz_1$.  Then we may write 
\begin{eqnarray*}
\alpha&=&\alpha_0+\frac{\alpha_1}{z_1}+\frac{\alpha_2}{z_1^2}+\cdots+\frac{\alpha_{p+1}}{z_{1}^{p+1}} \\
\beta&=&\beta_0+\frac{\beta_1}{z_1}+\frac{\beta_2}{z_1^2}+\cdots+\frac{\beta_{p+1}}{z_1^{p+1}},
\end{eqnarray*}

\noindent
where $\alpha_i,\beta_i$\,\,\,($i\geq 1$) do not involve $z_1$ and $dz_1$, and $\alpha_i,\beta_i$\,\,\, ($i\geq 0$) have poles of arbitrary order (at most) along $Y^\prime$ as their only singularities.  Since $d\varphi=0$, we have
\begin{equation*}
\begin{split}
d\varphi&= 
-dz_1\wedge d\alpha_0+d\beta_0- \frac{dz_1\wedge d\alpha_1+d\beta_1}{z_1}
-\frac{dz_1\wedge(d\alpha_2+\beta_1)-d\beta_2}{z_1^2} \\
&\text{\hskip 50pt} -\cdots-\frac{dz_1\wedge(d\alpha_{p+1}+p\beta_{p})-d\beta_{p+1}}{z_1^{p+1}} -(p+1)\frac{dz_1\wedge \beta_{p+1}}{z_1^{p+2}} \\
&=0.
\end{split}
\end{equation*}

\noindent
Hence, 
\begin{eqnarray}\label{e4.4}
& & d\alpha_1=d\alpha_2+\beta_1=d\alpha_3+2\beta_2=\cdots=d\alpha_{p+1}+p\beta_p=0,\\
& & (p+1)\beta_{p+1}=0, \nonumber \\
& & d\beta_1= d\beta_2=\cdots=d\beta_{p+1}=0, \nonumber\\ 
& & d\varphi_0=0, \quad \text{where}\quad \varphi_0=dz_1\wedge \alpha_0+\beta_0.\nonumber\end{eqnarray}

\noindent
Put 
\begin{equation*}
\theta=-\frac{\alpha_2}{z_1}-\frac{\alpha_3}{2z_1^2}-\cdots-\frac{\alpha_{p+1}}{pz_1^p},
\end{equation*}

\noindent
then 
\begin{equation}\label{e4.5}
\varphi=d\theta+\frac{dz_1}{z_1}\wedge \alpha_1+\varphi_0, \quad \text{and}
\quad d\varphi_0=0.
\end{equation}
\smallskip

\noindent
Hence if $p\geq 3$, since $d\alpha_1=d\varphi_0=0$, there exist local cross-sections $\gamma$ of $\Omega^{p-2}(\ast Y^\prime)$ and $\varphi_1$ of $\Omega^{p-1}(\ast Y^\prime)$ with $d\gamma=\alpha_1$ and $d\varphi_1=\varphi_0$ in a neighborhood of $x_0$. Put 
\begin{equation*}
\theta_1=\frac{dz_1}{z_1}\wedge \gamma+\varphi_1,
\end{equation*}

\noindent
then $\theta+\theta_1$ is a local cross-section of $\Omega^{p-1}(pY+\ast Y^\prime)$ and $\varphi=d(\theta+\theta_1)$.  This shows that 
$\mathcal{H}^p(\Omega_X^\cdot((1+\cdot)Y+\ast Y^\prime))=0$ \quad for \quad $p\geq 3$.  
If $p=2$, $\alpha_1$ in the expression (\ref{e4.5}) of $\varphi$ is a local cross-section of $\Phi^1(\ast Y^\prime)$.  Hence, as shown in the proof of Lemma \ref{lem2.1} (ii)-(c), there exists a constant $\lambda^\in \mathbb{C}$ and a local cross-section $\gamma$ of $\Omega_X^0(\ast Y^\prime)$ with 
\begin{equation*}
\alpha_1=\lambda\frac{dz_2}{z_2}+d\gamma.
\end{equation*}

\noindent
Furthermore, since $\varphi_0$ is a localcross-section of 
$\Phi^2(\ast Y^\prime)$, by Lemma \ref{lem2.1} (i)-(a), there exists a local cross-section $\varphi_1$ of $\Omega^1(\ast Y^\prime)$ with $d\varphi_1=\varphi_0$.  Put 
\begin{equation*}
\theta_1=\frac{dz_1}{z_1}\wedge \gamma+\varphi_1
\end{equation*}

\noindent
then $\theta+\theta_1$ is a local cross-section of $\Omega^1(2Y+\ast Y^\prime)$ at $x_0$ and 
\begin{eqnarray*}
\varphi&=& d\theta+\lambda\frac{dz_1\wedge dz_2}{z_1z_2}+\frac{dz_1}{z_1}\wedge d\gamma+\varphi_0 \\
 &=& \lambda\frac{dz_1\wedge dz_2}{z_1z_2}+d(\theta+\theta_1)
\end{eqnarray*}

\noindent
This shows that
\begin{equation*}
 \mathcal{H}^2(\Omega_X^\cdot((1+\cdot)Y+\ast Y^\prime))_{x_0}\simeq 
\mathbb{C}\{\frac{dz_1\wedge dz_2}{z_1z_2}\}. 
\end{equation*}

\noindent
If $p_1=1$, $\alpha_1$ is a meromorphic function, hence $d\alpha_1=0$ implies that $\alpha_1=\lambda$, a constant.  Since $\varphi_0$ is a local cross-section of $\Phi^1(\ast Y^\prime)$, by Lemma \ref{lem2.1} (ii)-(c), there exists $\varphi_1\in \Omega^0(\ast Y^\prime)_{x_0}$ such that 
\begin{equation*}
\varphi_0=\mu\frac{dz_2}{z_2}+d\varphi_1.
\end{equation*}
\smallskip

\noindent
Hence the expression of $\varphi$ in (\ref{e4.5}) becomes 
\begin{equation*}
\varphi=\lambda \frac{dz_1}{z_1}+\mu\frac{dz_2}{z_2}+d(\varphi_1+\theta)
.
\end{equation*}

\noindent
Since $\varphi_1+\theta\in \Omega^1(Y+\ast Y^\prime)$, this shows 
\begin{equation*}
\mathcal{H}^1(\Omega_X^\cdot((1+\cdot)Y+\ast Y^\prime))\simeq \mathbb{C}\{\frac{dz_1}{z_1},\frac{dz_2}{z_2}\}.
\end{equation*}
\smallskip

\noindent
$\mathcal{H}^0(\Omega_X^\cdot((1+\cdot)Y+\ast Y^\prime))\simeq \mathbb{C}_X$ is obvious.  To prove the same for $\mathcal{H}^p(\Omega_X^\cdot(\log Y+(1+\cdot)Y^\prime))$ is rather easy.  If $\varphi$ is a local cross-section of $\Phi^p(\log Y+(p+1)Y^\prime)$ in a neighborhood of $x_0$, then 
$\varphi$ is written as 
\begin{equation*}
\varphi=\frac{dz_1}{z_1}\wedge \alpha+\beta,
\end{equation*}

\noindent
where $\alpha\in \Omega^{p-1}((p+1)Y^\prime)$, $\beta\in \Omega^{p}((p+1)Y^\prime)$ do not involve $dz_1$.  Furthermore, we may assume that $\alpha$ does not 
involve $z_1$.  Then $d\varphi=$ implies $d\alpha=d\beta=0$, and by the same arguments as in the case of $\Omega_X^\cdot((p+1)Y+\ast Y^\prime)$, we can show that (\ref{e4.3}) for $\mathcal{H}^p(\Omega_X^\cdot(\log Y+(p+1)Y))$.
\end{proof}
\smallskip

\begin{lem}\label{lem4.3}
Assume we are under the setting at the begining of this section.  Particularly, we assume that the following conditions are satisfied:
\begin{eqnarray*}
H^p(X,\Omega_X^p(kY))&=&0,\\
H^p(Y,\Omega_Y^p(kZ))&=&0 \quad \text{for} \quad p\geq 1, q\geq 0, k\geq 1.
\end{eqnarray*}

\noindent
Then we have 
\begin{equation*}
H^p(X,\Omega_X^q(\log Y+(q+1)Y^\prime))=0 \quad \text{for} \quad p\geq 1, q\geq 0.
\end{equation*}
\end{lem}

\begin{proof} We consider the following exact sequence
\begin{equation*}
0\to \Omega_X^q\to \Omega_X^q(\log Y)\xrightarrow{R} \Omega_Y^{q-1}\to 0 
\quad (q\geq 1),
\end{equation*}

\noindent
where $R$ is the r{\'e}sidue map (cf. Lemma \ref{lem2.1} (ii)-(c)).  Tensoring 
$\mathcal{O}_X((q+1)Y^\prime)$ to this exact sequence, we have 
\begin{equation*}
0\to \Omega_X^q((q+1)Y^\prime)\to \Omega_X^q(\log Y+(q+1)Y^\prime)\to \Omega_Y^{q-1}((q+1)Z)\to 0.
\end{equation*}

\noindent
From the long exact sequence of cohomology associated to this sequence, the assertion of the lemma follows.
\end{proof}
\bigskip 
 We define 
\begin{eqnarray*}
I^p(X,\log Y+(p+1)Y^\prime)&:=&\frac{\Gamma(X,\Phi_X^p(\log Y+(p+1)Y^\prime))}{d\Gamma(X,\Omega_X^{p-1}(\log Y+pY^\prime))}, \\
I^p(X,(p+1)Y+\ast Y^\prime)&:=&\frac{\Gamma(X,\Phi_X^p((p+1)Y+\ast Y^\prime))}{d\Gamma(X,\Omega_X^{p-1}(pY+\ast Y^\prime))}.
\end{eqnarray*}

\noindent
Combining Proposition \ref{prop4.2} with Lemma \ref{lem4.3} implies the following:
\smallskip

\begin{prop}\label{prop4.4}
Assume that we are under the setting at the bigining of this section. Then 
\begin{equation*}
I^p(X,\log Y+(p+1)Y^\prime)\simeq I^p(X,(p+1)Y+\ast Y^\prime) \quad \text{for} \quad p\geq 0.
\end{equation*}
\end{prop}
\bigskip

We are now ready to define the R{\'e}sidue map 
\begin{equation*}
R\acute{e}s: I^p(X,(p+1)Y)\to I^{p-1}(Y,pZ)_0
\end{equation*}

\noindent
Let $\omega\in \Gamma(X,\Phi_X^p((p+1)Y))$ be given.  We think of $\omega$ as 
an element of $\Gamma(X,\Phi_X^p((p+1)Y+\ast Y^\prime)$. Then, by Propostion \ref{prop4.4}, there exists a $\varphi\in \Gamma(X,\Omega_X^{p-1}(pY+\ast Y^\prime)))$ such that $\omega-d\varphi\in\Gamma(\Phi_X^p(\log Y+(p+1)Y^\prime))$.  We take an 
open covering $\{U_i\}_{i\in I}$ of $X$ such that there is a local coordinate 
system $(z_1^i,\cdots,z_{n+1}^i)$ on each $U_i$, satisfying the following conditions:
\begin{equation}\label{e4.6}
\begin{split}
 {\rm (a)}&\,\,\,\text{If} \quad U_i\cap Y\not=\emptyset, z_1^i=0 \quad \text{is a defining equation of} \hskip 5pt  Y \hskip 5pt \text{in}\hskip 5pt U_i.\\ 
 {\rm (b)}&\,\,\, \text{If}\quad U_i\cap (Y\cap Y^\prime)\not=\emptyset, z_1^i=0\quad \text{and}\quad z_2^i=0 \quad \text{are defining equations}\\
& \text{of $Y$ and $Y^\prime$ in $U_i$, respectively}.
\end{split}
\end{equation}

\noindent
In each $U_i$ with $U_i\cap Y\not=\emptyset$, we can write $\omega-d\varphi$ as 
\begin{equation}\label{e4.7}
\omega-d\varphi=\frac{dz_1^i}{z_1^i}\wedge \alpha_i+\beta_i,
\end{equation}

\noindent
where $\alpha_i\in \Gamma(U_i,\Phi_X^{p-1}((p+1)Y^\prime))$, $\beta_i\in \Gamma(U_i,\Phi_X^p((p+1)Y^\prime))$, $\alpha_i$ and $\beta_i$ does not involve $dz_1^i$.  We can easily see $\alpha_{i\vert Y}=\alpha_{j\vert Y}$ if $U_i\cap U_j\cap Y\not=\emptyset$, hence $\{\alpha_{i\vert Y}\}$ defines an element of $\Gamma(Y,\Phi_X^{p-1}((p+1)Z))$.@We claim that $\{2\pi \sqrt{-1}\alpha_{i\vert Y}\}$ determine a unique element of $I^{p-1}(Y,(p+1)Z))$, not depending on the chice of $\varphi$.  In fact, if $\varphi^\prime$ is another element of $\Gamma(X,\Omega_X^{p-1}((p+1)(Y+Y^\prime))$ with 
$\omega-d\varphi^\prime\in \Gamma(X,\Phi_X^p(\log Y+(p+1)Y^\prime))$ and 
\begin{equation*}
\omega-d\varphi^\prime=\frac{dz_1^i}{z_1^i}\wedge \alpha_i^\prime+
\beta_i^\prime
\end{equation*}

\noindent
is the expression of $\omega-d\varphi^\prime$ as in (\ref{e4.7}), then 
\begin{equation*}
d(\varphi^\prime-\varphi)=\frac{dz_1^i}{z_1^i}\wedge(\alpha_i-
\alpha_i^\prime)+(\beta_i-\beta_i^\prime)\in \Gamma(X,\Phi_X^p(\log Y+(p+1)Y^\prime)
\end{equation*}

\noindent
is zero in $I^p(X,\log Y+(p+1)Y^\prime)$. Hence, by Proposition \ref{prop4.4}, there exists an element $\psi\in \Omega_X^{p-1}(\log Y+pY^\prime))$ such that 
$d\psi=d(\varphi^\prime-\varphi)$.  Let 
\begin{equation*}
\psi=\frac{dz_1^i}{z_1^i}\wedge \gamma_i+\delta_i
\end{equation*}

\noindent
be the expression of $\psi$ as in (\ref{e4.7}).  Then, since $d\psi=d(\varphi^\prime-\varphi)$, we have 
\begin{equation}\label{e4.8}
d\gamma_{i\vert Y}=d_Y(\gamma_{i\vert Y})=\alpha_{i\vert Y}-\alpha_{i\vert Y}^\prime
\end{equation}

\noindent
for each $i$ with $U_i\cap Y\not=\emptyset$, where $d_Y$ denotes the exterior derivative on $Y$. Since $\{\gamma_{i\vert Y}\}$ is a global cross-section of 
$\Gamma(Y,\Omega_X^{p-1}(pZ))$, (\ref{e4.8}) shows that $\{\alpha_{i\vert Y}\}=\{\alpha^\prime_{i\vert Y}\}$ in $I^{p-1}(Y,(p+1)Z)$.  Furthermore, the arguments above also show that if $\omega$ is a derived form, then so is $\{\alpha^\prime_{i\vert Y}\}$.  Therefore, we conclude that the correspondence 
\begin{equation*}
\omega \longmapsto \{\alpha^\prime_{i\vert Y}\}
\end{equation*}

\noindent
determine a map $I^{p}(X,(p+1)Y)\to I^{p-1}(Y,(p+1)Z)$.  Since $I^{p-1}(Y,(p+1)Z)\simeq I^{p-1}(Y,pZ)$ by Proposition \ref{prop2.3}, this map is thought of 
as a map from $I^{p}(X,(p+1)Y)$ to $I^{p-1}(Y,pZ)$, which we define to be the {\it generalized Poncar{\'e}\/} r{\'e}sidue map and denote it {\it R{\'e}s\/}.  We denote $\{\alpha^\prime_{i\vert Y}\}$ by $r\acute{e}s[\omega]$ (determined up to derived forms) and call {\it r{\'e}sidue form\/} of $\omega$. 
\smallskip
\begin{prop}\label{prop4.5}
\begin{equation*}
R\acute{e}s(I^p(X,(p+1)Y))\subset I^{p-1}(Y,pZ)_0
\end{equation*}
\end{prop}

\begin{proof}
For a $\omega\in \Gamma(X,\Phi^p((p+1)Y))$, we shall show that its r{\'e}sidue 
form $r\acute{e}s[\omega]=\{\alpha_{i\vert Y}\}$ (precisely speaking, a closed form
representing the class $r\acute{e}s[\omega]$ of $I^{p-1}(Y,(p+1)Z))$ is of the second kind in the sense of Picard-Lefshetz.  From this the assertion of the proposition follows, since $I^{p-1}(Y,(p+1)Z)_0\simeq I^{p-1}(Y,pZ)_0$.  As before we take an open covering $\{U_i\}_{i\in I}$ of $X$ such that there is a local coordinate 
system $(z_1^i,\cdots,z_{n+1}^i)$ on each $U_i$, subject to the conditions in 
(\ref{e4.6}), and take a $\varphi\in \Gamma(X,\Omega^{p-1}(pY+\ast Y^\prime))$ such that $\omega-d\varphi\in \Gamma(X,\Phi_X^{p}(\log Y+(p+1)Y^\prime))$.  On each
$U_i$ with $U_i\cap Y\not=\emptyset$, we write 
\begin{equation}\label{e4.9}
\omega-d\varphi=\frac{dz_1^i}{z_1^i}\wedge \alpha_i+\beta_i
\end{equation}

\noindent
as in (\ref{e4.7}).  We will show that for a point $x_0\in Z\cap U_i$, 
$r\acute{e}s[\omega]_{\vert U_i}=\alpha_{i\vert Y}$ is a holomorphic form 
modulo derived 
meromorphic forms in a sufficiently small neighborhood of $x_0$ in $Y$.  For 
this end we take a generic prime hypersurface section $Y^{\prime\prime}$ which 
is linearly equivalent to $Y^\prime$, which does not go through $x_0$ and intersect 
$Y$ and $Y^\prime$ transversely.  We think $\omega$ as an element of 
$\Gamma(X,\Phi^p((p+1)Y+\ast Y^{\prime\prime})))$.  Since $I^p(X,(p+1)Y+\ast Y^{\prime\prime}))\simeq I^p(X,\log Y+(p+1)Y^{\prime\prime})$ by Proposition \ref{prop4.4}, there exists a $\varphi^\prime\in \Gamma(X,\Omega^{p-1}(pY+\ast Y^{\prime\prime}))$ with $\omega-d\varphi^\prime\in \Gamma(X,\Phi^p(\log Y+(p+1)Y^{\prime\prime}))$.  Let 
\begin{equation}\label{e4.10}
\omega-d\varphi^\prime=\frac{dz_1^i}{z_1^i}\wedge \alpha_i^\prime+\beta_i^\prime
\end{equation}

\noindent
be the expression of $\omega-d\varphi^\prime$ as in (\ref{e4.7}) on each 
$U_i\cap Y\not=\emptyset$.  If $U_{i_0}$ is the coordinate neighborhood with 
$x_0\in U_{i_0}\cap Z$, since $Y^{\prime\prime}$ does not go through $x_0$, 
$\alpha_{i_0\vert Y}^{\prime}$ is holomorphic in a sufficiently open neighborhood of $x_0$ in $U_{i_0}\cap Y$.  From (\ref{e4.9}) and (\ref{e4.10}), 
\begin{equation}\label{e4.11}
d(\varphi^\prime-\varphi)=\frac{dz_1^{i_0}}{z_1^{i_0}}\wedge (\alpha_{i_0}-\alpha_{i_0}^\prime)+(\beta_{i_0}-\beta_{i_0}^\prime).
\end{equation}

\noindent
Since $d(\varphi^\prime-\varphi)\in \Gamma(X,\Phi^p(\log Y+\ast(Y^\prime+Y^{\prime\prime})))$ is zero in $I^p(X,(p+1)Y+\ast (Y^\prime+Y^{\prime\prime}))$, by Proposition \ref{e4.4}, there exists a $\psi\in \Gamma(X,\Omega_X^{p-1}(\log Y+\ast(Y^\prime+Y^{\prime\prime})))$ with $d\psi=d(\varphi^\prime-\varphi)$. On each $U_i$, we write 
\begin{equation}\label{e4.12}
\psi=\frac{dz_1^{i}}{z_1^{i}}\wedge \gamma_i+\xi_{i_0}
\end{equation}

\noindent
as in (\ref{e4.7}).@Then $d\psi=d(\varphi^\prime-\varphi)$ implies 
\begin{equation*}
d\gamma_i=\alpha_i-\alpha_i^\prime.
\end{equation*}

\noindent
Hence $d_Y(\gamma_{i\vert Y})=\alpha_{i\vert Y}-\alpha_{i\vert Y}^\prime$ for each $i$ where $d_Y$ denotes the exterior derivation on $Y$.  This means $d_Y(r\acute{e}s[\psi])=r\acute{e}s[\alpha]-r\acute{e}s[\alpha^\prime]$ where 
$r\acute{e}s[\psi]\in \Gamma(Y,\Omega_Y^{p-2}(p(Y^\prime+Y^{\prime\prime})))$.  Since 
$r\acute{e}s[\alpha^\prime]$ is holomorphic at $x_0$, so is $r\acute{e}s[\alpha]$  modulo derived meromorphic forms as requied.
\end{proof}
\smallskip

{\bf Proof of Theorem 4.1:} \par
\smallskip
 We can now easily deduce Theorem \ref{thm4.1} from 
what we have proved till now.  First, by Theorem \ref{thm1.15},
\begin{equation*}
H^n(Y,\mathbb{C})_0=R^{n+1}(H^{n+1}(X-Y,\mathbb{C}))\oplus r^n(H^n(X,\mathbb{C})_0).
\end{equation*}

\noindent
By Theorem \ref{thm3.1}, 
\begin{equation}\label{e4.13}
F^k(H^n(Y,\mathbb{C})_0=R^{n+1}(F^{k+1}H^{n+1}(X-Y,\mathbb{C}))\oplus r^n(F^kH^n(X,\mathbb{C})_0).
\end{equation}

\noindent
By Theorem \ref{thm3.3}, (\ref{e3.5}), 
\begin{equation}\label{e4.14}
F^{k+1}H^{n+1}(X-Y,\mathbb{C}))\simeq I_{k+1}^{n+1}(X,(n+2)Y).
\end{equation}

\noindent
Applying Theorem \ref{thm3.3}, (\ref{e3.6}) to the pair $(X,Y^\prime)$ instead of $(X,Y)$, we have 
\begin{equation}\label{e4.15}
F^{k}H^{n}(X,\mathbb{C}))_0\simeq I_{k}^{n}(X,(n+1)Y^\prime)_0.
\end{equation}

\noindent
From (\ref{e4.13}), (\ref{e4.14}) and (\ref{e4.15}) it follows that
\begin{equation*}
F^{k}H^{n}(Y,\mathbb{C}))_0=R^{n+1}I_{k+1}^{n+1}(X,(n+2)Y^\prime)\oplus r^n(I^n_k(X,(n+1)Y^\prime)_0).
\end{equation*}

\noindent
Here the map $R^{n+1}:I^{n+1}(X,(n+2)Y)\simeq H^{n+1}(X-Y,\mathbb{C})\to H^n(Y,\mathbb{C})$ should be interpreted in terms of $C^\infty$ De Rham group as follows: By use of isomorphisms
\begin{equation*}
H^{n+1}(X-Y,\mathbb{C})\simeq \mathbb{H}^{n+1}(X,\Omega_X^\cdot(\log Y))\simeq 
I^{n+1}(X,(n+2)Y)\simeq H^{n+1}(A^\cdot(\log Y)),
\end{equation*}

\noindent
(cf. Proposition \ref{prop2.2} and its proof), we can take a $\varphi\in {\rm Ker}\{(A^{n+1}(\log Y))\to A^{n+2}(\log Y)\}$ with $\omega=\varphi$ modulo 
$dA^n(\log Y)$ for a $\omega\in \Gamma(X,\Phi^{n+1}((n+2)Y))$.  $\varphi$ is written as
\begin{equation*}
\varphi=\alpha\wedge \eta+\beta,
\end{equation*}

\noindent
where $\eta$ is $C^\infty$ form of type $(1,0)$ with the property $\overline{\partial}\eta$ represents the first Chern class $c_1([Y])$, and $\alpha\in A^{n-1}(X)$, $\beta\in A^{n+1}(X)$ (cf. (\ref{e1.4}). $d\varphi=0$ implies $d_Y(\alpha_{\vert Y})=0$.  Then $R^{n+1}([\omega])$ ($[\omega]\in I^{n+1}(X,(n+2)Y))$ is defined by 
\begin{equation*}
R^{n+1}([\omega])=2\pi\sqrt{-1}\Big[\alpha_{\vert Y}\Big],
\end{equation*}

\noindent
where $[\alpha_{\vert Y}]$ denote the De Rham cohomology class represented by 
$\alpha_{\vert Y}$.  Taking into consideration this fact, we will be done 
if we see 
\begin{equation}\label{e4.16}
R^{n+1}(I^{n+1}_{k+1}(X,(n+2)Y))=R\acute{e}s(I^{n+1}_{k+1}(X,(n+2)Y)
\end{equation}

\noindent
in the De Rham cohomology. To see  this, we first note that both of the right and left hand sides of (\ref{e4.16}) are included in $H^n(Y,\mathbb{C})_0$. due to Theorem \ref{thm1.15} and Theorem \ref{thm2.4}.  Hence, by Proposition \ref{prop1.9} and Proposition \ref{prop1.10}, in order to prove (\ref{e4.16}), it suffices to show that 
\begin{equation}\label{e4.17}
\int_{\tau_{\varepsilon}(\gamma)}\omega=\int_{\gamma}r\acute{e}s[\omega]
\end{equation}

\noindent
for a $\omega\in \Gamma(X,\Phi^{n+1}((n+2)Y))$ and an $n$ cycle $\gamma$ lying in $Y-Z$, where $\tau_{\varepsilon}(\gamma)$ is $\partial U_{\varepsilon\vert \gamma}$, the restriction of the boundary of a topological {\it $\varepsilon$ tublorneighborhood\/} $U_{\varepsilon}$ of $Y$ in $X$ to $\gamma$.  We are now going to prove (\ref{e4.17}). We take the local expression (\ref{e4.7}) of $\omega$ with respect to some open covering $\{U_i\}_{i\in I}$ of $X$ and a local coordinate system $(z_1^i,\cdots,z_{n+1}^i)$ on each $U_i$, subject to the conditions in 
(\ref{e4.6}).  Let $\{\rho_i\}$ be a partition of unity subordinate to the covering $\{U_i\}_{i\in I}$.  Then 
\begin{eqnarray*}
\int_{\tau_{\varepsilon}(\gamma)}\omega&=& 
\int_{\tau_{\varepsilon}(\gamma)}\sum_{i}\rho_i(\frac{dz_1^i}{z_1^i}\wedge \alpha_i+\beta_i)+d\varphi \\
&=& \int_{\tau_{\varepsilon}(\gamma)}\sum_{i}\rho_i(\frac{dz_1^i}{z_1^i}\wedge \alpha_i+\beta_i) \\
&=& \sum_i\int_{\tau_{\varepsilon}(\gamma)}\rho_i(\frac{dz_1^i}{z_1^{i}}\wedge \alpha_i+\beta_i).
\end{eqnarray*}

\noindent
Locally, $\tau_{\varepsilon}(\gamma)$ looks like $\mathbb{R}^{n+1}\times \{~\vert z\vert=\varepsilon~\vert~z\in \mathbb{C}~\}$ ($\varepsilon>0$).  Hence 
\begin{eqnarray*}
 \sum_i\int_{\tau_{\varepsilon}(\gamma)\cap U_i}\rho_i(\frac{dz_1^{(i)}}{z_1^{(i)}}\wedge \alpha_i+\beta_i)
 &=& \lim_{\varepsilon\to 0}\sum_i\int_{\tau_{\varepsilon}(\gamma)\cap U_i}\rho_i(\frac{dz_1^{i}}{z_1^{i}}\wedge \alpha_i+\beta_i)\\
 &=& 2\pi\sqrt{-1}\sum_i(\rho_i\alpha_i)_{\vert\gamma\cap U_i} \\
 &=& 2\pi\sqrt{-1}r\acute{e}s[\omega]
 \end{eqnarray*}
 
 \noindent
 as required.  This completes the proof of Theorem \ref{thm4.1}.
 \smallskip
 
 \begin{rem}\label{rem4.1}
 For $[\omega]\in I_{k+1}^{n+1}(X,(n+2)Y)$ it can be proved more directly that the Hodge type of $R^{n+1}([\omega])=R\acute{e}s([\omega])$ is $(n,0)+(n-1,1)+\cdots+(k,n-k)$.  By virtue of the isomorphism 
\begin{eqnarray*}
& & I_{k+1}^{n+1}(X,(n+2)Y)\simeq H^{n+1}(F^{\prime k+1}(B^\cdot)) \\
& & \\
&=& \frac{{\rm Ker}\{\sum_{\ell=0}^{n-k}B_X^{n-\ell+1,\ell}(n-\ell-k+1)\xrightarrow{d}  
\sum_{\ell=0}^{n-k+1}B_X^{n-\ell+2,\ell}(n-\ell-k+2)\}}
{{\rm Im}\{\sum_{\ell=0}^{n-k-1}B_X^{n-\ell,\ell}(n-\ell-k)\xrightarrow{d}
\sum_{\ell=0}^{n-k}B_X^{n-\ell+1,\ell}(n-\ell-k+1)\}}
\end{eqnarray*}

\noindent
(cf. the proof of Theorem \ref{thm3.3}, (\ref{e3.5})), $\omega\in \Gamma(X,\Phi_X^{n+1}((n+2)Y)$ is cohomologous to a closd form $\varphi$ of $\sum_{\ell=0}^{n-k}B_X^{n-\ell+1,\ell}(n-\ell-k+1)$ in the De Rham cohomology.  If we wtite $\varphi$ as 
\begin{equation*}
\varphi=\varphi^{(n+1,0)}+\varphi^{(n,1)}+\cdots+\varphi^{(k+1,n-k)},
\end{equation*}

\noindent
where $\varphi^{(n-\ell+1,\ell)}\in B_X^{n+\ell-1,\ell}(n-\ell-k+2)$ \quad 
($0\leseq \ell\leseq n-k$), then each $\varphi^{(n-\ell+1,\ell)}$  
is written in each $U_i$ as 
\begin{equation*}
\varphi^{(n-\ell+1,\ell)}=\frac{\alpha_i^{(n-\ell,\ell)}dz_1^{i}}{(z_1^{i})^{n-\ell-k+1}}+\frac{\beta_i^{(n-\ell+1,\ell)}}{(z_1^{i})^{n-\ell-k}}
\end{equation*}

\noindent
where $\alpha_i^{(n-\ell,\ell)}$, $\beta_i^{(n-\ell+1,\ell)}$ are regular $C^\infty$ differential forms of types $(n-\ell,\ell)$, $(n-\ell+1,\ell)$, respectively, not involving $z_1^i$, where $z_1^i=0$ is the local defining equation of $Y$. This is because $(z_1^i)^{n-\ell-k}\varphi^{(n-\ell+1,\ell)}$ and $(z_1^i)^{n-\ell-k}dz_1^i\wedge\varphi^{(n-\ell+1,\ell)}$ are $C^\infty$ regular forms 
by the definition of $B_X^{n-\ell+1,\ell}(n-\ell-k+1)$. Put
\begin{equation*}
\psi_i^{(n-\ell,\ell)}:=\frac{\alpha_i^{(n-\ell,\ell)}}{(n-\ell-k)(z_1^i)^{n-\ell-k}} \qquad (0\leseq \ell \leseq n-k-1),
\end{equation*}

\noindent
then 
\begin{eqnarray*}
\eta_i^{(n-\ell+1,\ell)+(n-\ell,\ell+1)}&:=&d\psi_i^{(n-\ell,\ell)}+\varphi^{(n-\ell+1,\ell)} \\
&=& \frac{d\alpha_i^{(n-\ell,\ell)}}{(n-\ell-k)(z_1^i)^{n-\ell-k}}+\frac{\beta_i^{(n-\ell+1,\ell)}}{(z_1^i)^{n-\ell-k}}
\end{eqnarray*}

\noindent
is a semi-meromorphic form of type $(n-\ell+1,\ell)+(n-\ell,\ell+1)$ and has poles of order $n-\ell-k$ along $Y$. Let $\{\rho_1\}$ be a partition of unity subodinate to the open covering $\{U_i\}_{i\in I}$ as before.  We put 
\begin{eqnarray*}
\psi^{(n-\ell,\ell)}&=& \sum_i\rho_i\psi_i^{(n-\ell,\ell)}, \\
 & & \\
\eta^{(n-\ell+1,\ell)+(n-\ell,\ell+1)}&=& \sum_i\rho_i\eta_i^{(n-\ell+1,\ell)+(n-\ell,\ell+1)}
\end{eqnarray*}

\noindent
Now, 
\begin{eqnarray*}
\varphi^{(n-\ell+1,\ell)}-d\psi^{(n-\ell,\ell)}&=& 
\varphi^{(n-\ell+1,\ell)}-\sum_id\rho_i\psi_i^{(n-\ell,\ell)}+\sum_i\rho_id\psi_i^{(n-\ell,\ell)}\\
&=& \sum_i\rho_i\eta_i^{(n-\ell+1,\ell)+(n-\ell,\ell+1)}-\sum_id\rho_i\psi_i^{(n-\ell,\ell)}\\
&=& \eta^{(n-\ell+1,\ell)+(n-\ell,\ell+1)}-\sum_id\rho_i\psi_i^{(n-\ell,\ell)}
\end{eqnarray*}

\noindent
which is a semi-morphic form of type $(n-\ell+1,\ell)+(n-\ell,\ell+1)$ having poles of order $n-\ell-k$ along $Y$.  Continuing this process, $\varphi^{(n-\ell+1,\ell)}$ \quad ($0\leseq \ell\leseq n-k$) is reduced to a semi-meromorphic form of thpe $(n-\ell+1,\ell)+\cdots+(k+1,n-k)$, having poles of order $1$ along $Y$ modulo derived forms.  Hence $\varphi$ is reduced to a closed semi-meromorphic 
form $\xi$ of $A^{n+1,0}(\log Y)+\cdots+A^{k+1,n-k}(\log Y)$ modulo derived forms. Hence the Hodge type of $R^{n+1}([\omega])=R^{n+1}([\xi])$ is $(n,0)+(n-1,1)+\cdots+(k,n-k)$. 
\end{rem}


\end{document}